\documentclass[a4paper,12pt]{article}
\usepackage{amsmath, cases}
\usepackage{amscd}
\usepackage[xdvi]{graphicx}
\usepackage[dvips]{color}
\usepackage{txfonts}

\newcommand{\bm}[1]{\mbox{\boldmath $#1$}}

\makeatletter

\@addtoreset{equation}{section}

\newcounter{def}[section]
\renewcommand{\thedef}{\stepcounter{def}\thesection.\@arabic\c@def }

\begin{document}
\setlength{\baselineskip}{24pt}
\begin{center}
\textbf{\LARGE{Reduction of weakly nonlinear parabolic partial differential equations}}
\end{center}

\setlength{\baselineskip}{14pt}

\begin{center}
Institute of Mathematics for Industry, Kyushu University, Fukuoka,
819-0395, Japan

\large{Hayato CHIBA} \footnote{E mail address : chiba@imi.kyushu-u.ac.jp}
\end{center}
\begin{center}
Apr 5, 2013
\end{center}

\begin{center}
\textbf{Abstract}
\end{center}

It is known that the Swift-Hohenberg equation
$\partial u/\partial t = -(\partial _x^2 + 1)^2u + \varepsilon (u-u^3)$ can be reduced 
to the Ginzburg-Landau equation (amplitude equation) $\partial A/\partial t = 4\partial _x^2 A + \varepsilon (A-3A|A|^2)$
by means of the singular perturbation method.
This means that if $\varepsilon >0$ is sufficiently small,
a solution of the latter equation provides an approximate solution of the former one.
In this paper, a reduction of a certain class of a system of
nonlinear parabolic equations $\partial u/\partial t = \mathcal{P}u + \varepsilon f(u)$ is proposed.
An amplitude equation of the system is defined and an error estimate of solutions is given.
Further, it is proved under certain assumptions that if the amplitude equation has a stable steady state,
then a given equation has a stable periodic solution .
In particular, near the periodic solution, the error estimate of solutions holds uniformly in $t>0$.

\noindent \textbf{Keywords}: amplitude equation; renormalization group method; reaction diffusion equation

% \tableofcontents

%%%%%%%%%%%%%%%%%%%%%%%%%%%%%%%%%%%%%%%%%%%%%%%%%%%%%%%%%%%%%%%%%%%%%%%%%%%%%%%%%%%%%%%%%%%%%%%%%%
%%%%%%%%%%%%%%%%%%%%%%%%%%%%%%%%%%%%%%%%%%%%%%%%%%%%%%%%%%%%%%%%%%%%%%%%%%%%%%%%%%%%%%%%%%%%%%%%%%

\section{Introduction}

A reduction of a certain class of nonlinear parabolic partial differential equations (PDEs)
\begin{equation}
\frac{\partial u}{\partial t} = \mathcal{P}u + \varepsilon f(u), 
\quad u = u(t,x) \in \mathbf{C}^m, (t,x) \in \mathbf{R} \times \mathbf{R}^d,
\label{1-1}
\end{equation}
is considered, where $\varepsilon >0$ is a small parameter, 
$\mathcal{P}$ is an elliptic differential operator with constant coefficient and 
$f$ is a function on $\mathbf{C}^m$ satisfying suitable assumptions.
Our study is motivated by the following three problems.

\textbf{Case 1.} 
It is well known that the Swift-Hohenberg equation
\begin{equation}
\frac{\partial u}{\partial t} = -(\partial _x^2 + k^2)^2u + \varepsilon (u-u^3),
\quad u,x\in \mathbf{R},
\label{1-2}
\end{equation}
with a parameter $k\in \mathbf{R}$, can be reduced to the Ginzburg-Landau equation (amplitude equation)
\begin{equation}
\frac{\partial A}{\partial t} = 4k^2 \partial _x^2 A + \varepsilon (A-3A|A|^2),
\label{1-3}
\end{equation}
by means of the multiscaling method \cite{Col1} or the renormalization group (RG) method \cite{Che}.
Let $v_0$ be a function in some function space and give initial conditions
\begin{eqnarray*}
u(0,x) = v_0(\sqrt{\varepsilon }x) e^{ikx} + \overline{v_0(\sqrt{\varepsilon }x)} e^{-ikx},
\quad A(0,x) = v_0(\sqrt{\varepsilon }x),
\end{eqnarray*}
for Eqs.(\ref{1-2}) and (\ref{1-3}), respectively.
In \cite{Col2}, it is proved that there exists a positive number $C$ such that 
solutions of the two initial value problems  satisfy
\begin{equation}
|| u(t,x) - (A(t,x)e^{ikx} + \overline{A(t,x)}e^{-ikx}) || \leq C\sqrt{\varepsilon },
\label{1-4}
\end{equation}
up to the time scale $t\sim O(1/\varepsilon )$ with a certain norm.
In this case, a fourth-order PDE is reduced to a second-order PDE.

\textbf{Case 2.} 
Let $\Omega = \mathbf{R} \times (0,l)$ be the strip region on $\mathbf{R}^2$.
Consider the boundary value problem of a system of reaction diffusion equations on $\Omega $
\begin{equation}
\left\{ \begin{array}{l}
\displaystyle \frac{\partial u}{\partial t}
 = d(\partial ^2_{x}u + \partial ^2_{y}u) + ku-v + \varepsilon (u-u^3),  \\[0.3cm]
\displaystyle \frac{\partial v}{\partial t}
 = \partial ^2_{x}v + \partial ^2_{y}v + u-v,  \\[0.3cm]
\displaystyle \frac{\partial u}{\partial y}\Bigl|_{y=0,l} = \frac{\partial v}{\partial y}\Bigl|_{y=0,l} = 0 , 
\end{array} \right.
\label{1-14}
\end{equation}
where $l, d$ and $k$ are positive constants.
This type problem was introduced by Chen, Ei and Lin 
\cite{Ei} to investigate a stripe pattern observed in the skin of angelfish.
Under certain assumptions on parameters so that the system undergoes Turing instability,
they formally derived an amplitude equation of the form
\begin{equation}
\frac{\partial A}{\partial t}
 = - \frac{2d^2}{(k+d)(1-d)} \frac{\partial ^4A}{\partial x^4} 
     + \frac{\varepsilon }{1-d}(A - 3A|A|^2),
\label{1-17b}
\end{equation}
without any mathematical justification.
It is remarkable that the amplitude equation is a fourth-order equation while a given system is a second-order equation
because of a certain degeneracy of the dispersion relation.
On the other hand, a system of equations becomes a single equation and the number of space variables is reduced.

\textbf{Case 3.} 
Let us consider a system of reaction diffusion equations on $\mathbf{R}$
\begin{equation}
\left\{ \begin{array}{l}
\displaystyle \frac{\partial u}{\partial t}
 = D \partial ^2_{x}u + v + \varepsilon (u-u^3),  \\[0.3cm]
\displaystyle \frac{\partial v}{\partial t}
 = D \partial ^2_{x}v -u,  \\[0.3cm]
\end{array} \right.
\label{1-7b}
\end{equation}
where $D>0$ is a constant.
This system can be reduced to the Ginzburg-Landau equation
\begin{equation}
\frac{\partial A}{\partial t} = D \partial _x^2 A + \frac{\varepsilon }{2} (A-3A|A|^2).
\label{1-7c}
\end{equation}
In this case, the order of differential equations are the same, while a system is reduced to a single equation.

A purpose in this paper is to give a unified theory of such a reduction of PDEs (\ref{1-1}),
and give an error estimate of solutions of amplitude equations.
Furthermore, we will partially prove a conjecture by Oono and Shiwa \cite{Oono},
which states that if a given PDE is structurally stable, its amplitude equation provides
the qualitative features of the given system.
For example, a stable invariant manifold of an amplitude equation implies the existence of 
a stable invariant manifold of a given system.
For example, we will prove that Eq.(\ref{1-17b}) actually provides an approximate solution of the system (\ref{1-14}).
Further, a conjecture by \cite{Oono} is solved in the following sense;
if the amplitude equation (\ref{1-17b}) has a stable steady state, then the system (\ref{1-14})
has a corresponding stable periodic solution.

Since general results for (\ref{1-1}) is rather complicated,
we divide main results into several steps as follows:

\textbf{(1)} In this Introduction, our main results are stated for one-dimensional problems
$u\in \mathbf{C}$ and $x\in \mathbf{R}$ for simplicity.
\\[-0.2cm]

\textbf{(2)} In Sec.3, the asymptotic behavior of linear semigroups generated by elliptic differential operators are
investigated.

(2-i) Sec.3.1 deals with the case $u\in \mathbf{C}$ and $x\in \mathbf{R}^d$.
The asymptotic behavior of a semigroup $e^{\mathcal{P}t}$ is given 
under the assumptions (B1) to (B3) (Propositions 3.1, 3.2 and 3.3).

(2-ii) In Sec.3.2, the case $u\in \mathbf{C}^m$ and $x\in \mathbf{R}^d$ is considered.
The asymptotic behavior of a semigroup $e^{\mathcal{P}t}$ is given 
under the assumptions (C0) to (C3) (Proposition 3.6).
\\[-0.2cm]

\textbf{(3)} Sec.4 is devoted to nonlinear estimates and our main results are given.

(3-i) In Sec.4.1, the case $u\in \mathbf{C}$ and $x\in \mathbf{R}^d$ is considered,
which includes Case 1 above as an example.
The definition of an amplitude equation (reductive equation) is given.
An error estimate of solutions (Thm.4.2) is proved under the assumptions (D1) to (D3),
and the existence of stable periodic solutions (Thm.4.3) is proved under the assumptions (D1) to (D4).

(3-ii) In Sec.4.2, the case $u\in \mathbf{C}^m$ and $x\in \mathbf{R}^d$ is considered,
which includes Case 2 and Case 3 above.
An error estimate of solutions (Thm.4.13) is proved under the assumptions (E0) to (E3),
and the existence of stable periodic solutions (Thm.4.14) is proved under the assumptions (E0) to (E4).
Theorems 4.13 and 4.14 include all previous results.
\\[-0.2cm]

Since we need several integers to state our final results in Sec.4, we summarize some of them for 
the convenience of the reader.
\begin{itemize}
 \item An integer $m$ denotes the dimension of unknown function: $u \in \mathbf{C}^m$.
 \item An integer $d$ denotes the dimension of space variables:  $x \in \mathbf{R}^d$.
 \item An integer $M$ denotes the degeneracy of the dispersion relation, which determines the order of differentiation
of the amplitude equation. For Case 1 and 3, $M=2$, while $M=4$ for Case 2.
 \item An integer $D \, (1\leq D \leq d)$ denotes a dimension of the critical direction (see Sec.3.1 for the detail),
which gives the number of space variables included in the amplitude equation.
For Case 2, $d = 2$ and $D=1$. 
 \item An integer $N$ gives the number of critical wave numbers, at which the spectrum of the operator $\mathcal{P}$
is tangent to the imaginary axis.
In other words, it gives the dimension of a center subspace (see Sec.4.1 for the detail).
An amplitude equation becomes a system of $N$-equations.
\end{itemize}

Although our purpose is \textit{a system} of PDEs including a degenerate case as above, it would be better to start with
a one-dimensional case for the sake of simplicity.
In this Introduction, we suppose $u\in \mathbf{C}$ and $x\in \mathbf{R}$ are one-dimensional variables
to state our main results as simple as possible.
Higher dimensional problems will be treated after Sec.2.

Let $P(x) = \sum^q_{\alpha =0} a_\alpha x^\alpha $ be a polynomial of $x\in \mathbf{R}$ and 
$\mathcal{P}:= P(\partial _x)$ a corresponding differential operator on $\mathbf{R}$.
For the operator $\mathcal{P}$, we suppose the following:
\\[0.2cm]
\textbf{(A1)} $\mathrm{Re}[P(i\xi)] \leq 0$ for any $\xi \in \mathbf{R}$.
\\
\textbf{(A2)} There exist $\omega , k\in \mathbf{R} \, ((\omega ,k) \neq (0,0))$ and an integer $M$ such that 
\begin{eqnarray}
& & P(\pm ik) = \pm i\omega , \label{1-5}\\
& & P'(\pm ik) = \cdots  = P^{(M-1)}(\pm ik) = 0, \label{1-6}\\
& & P^{(M)}(ik) = P^{(M)}(-ik)\neq 0. \label{1-7}
\end{eqnarray}
\textbf{(A3)} $a_q i^q < 0$ and $P^{(M)} (ik) i^M < 0$.
\\

The assumption (A1) implies that the spectrum $\sigma (\mathcal{P})$ of $\mathcal{P}$ calculated 
in a suitable space, which coincides with $P(i\mathbf{R})$, is included in the closed left half plane.
If $\sigma (\mathcal{P})$ were included in the open left half plane, $u=0$ is linearly stable.
Since we are interested in a bifurcation occurred at $\varepsilon =0$, we supposed in Eq.(\ref{1-5})
that $\sigma (\mathcal{P})$ includes points $\pm i\omega $ on the imaginary axis.
The integer $M$ represents the degeneracy of the dispersion relation $\lambda = P(i\xi)$.
Define the operator $\mathcal{Q}$ to be
\begin{equation}
\mathcal{Q} = \frac{P^{(M)} (ik)}{M!} \frac{\partial ^M}{\partial x^M}.
\label{1-8}
\end{equation}
The assumption (A3) assures that $\mathcal{P}$ and $\mathcal{Q}$ are elliptic.
In this section, we further suppose that integers $j$ satisfying $P(ijk) = ij\omega$ are only $\pm 1$
(this will be removed after Sec.2).
For the Swift-Hohenberg equation, $\omega =0, M = 2$ and $k$ is the $k$ in Eq.(\ref{1-2}).

Let $B^r :=BC^r(\mathbf{R}; \mathbf{C})$ be a vector space
of complex-valued functions $f$ on $\mathbf{R}$ such that $f(x), f'(x),$ $\cdots ,f^{(r)}(x)$ are bounded uniformly continuous.
This is a Banach space with the norm defined by $|| f || = \sup \{ |f(x)|, \cdots , |f^{(r)}(x)|\}$.
For a function $f : B^r\to B^r$, define the function $R : (B^r \times B^r)\to B^r$ to be
\begin{equation}
R(A_1, A_2) = \left\{ \begin{array}{ll}
\displaystyle \frac{k}{2\pi} \int^{2\pi / k}_{0}\! f(A_1e^{ikx} + A_2 e^{-ikx}) e^{-ikx}dx & (\text{when}\,\, k\neq 0), \\
\displaystyle \frac{\omega }{2\pi} \int^{2\pi / \omega }_{0}\! 
      f(A_1e^{i\omega t} + A_2 e^{-i\omega t}) e^{-i\omega t}dt & (\text{when}\,\, \omega \neq 0).  \\
\end{array} \right.
\end{equation}
One can verify that these two expressions coincide with one another if $k\neq 0$ and $\omega \neq 0$.
For example if $f(u) = u-u^3$, then $R(A_1, A_2) = A_1 - 3A_1^2A_2$.

Now we consider two initial value problems
\begin{equation}
\displaystyle \frac{\partial u}{\partial t}  = \mathcal{P}u + \varepsilon f(u), \quad
\displaystyle u(0,x) = v_1(\eta x)e^{ikx} + v_2(\eta x)e^{-ikx} 
\label{1-10a}
\end{equation}
and
\begin{equation}
\left\{ \begin{array}{l}
\displaystyle \frac{\partial A_1}{\partial t} = \mathcal{Q}A_1 + \varepsilon R(A_1, A_2), \quad
\displaystyle A_1(0,x) = v_1(\eta x),\\[0.2cm]
\displaystyle \frac{\partial A_2}{\partial t} = \mathcal{Q}A_2 + \varepsilon R(A_2, A_1), \quad
\displaystyle A_2(0,x) = v_2(\eta x), \\
\end{array} \right.
\label{1-10b} 
\end{equation}
where $\eta := \varepsilon ^{1/M}$. We call the latter system the \textit{amplitude equation}.
When $\mathcal{P} = -(\partial _x^2 + k^2)^2$ and
$f(u) = u-u^3$, the Ginzburg-Landau equation (\ref{1-3}) is obtained as a special case $A_2 = \overline{A}_1$.
\\[0.2cm]
\textbf{Theorem \thedef.} Suppose $f : BC^r(\mathbf{R}; \mathbf{C}) \to BC^r(\mathbf{R}; \mathbf{C})\, (r\geq 1)$ is $C^1$ and $\varepsilon >0$ is sufficiently small.
For any $v_1, v_2 \in BC^r(\mathbf{R}; \mathbf{C})$, there exist positive numbers $C, T_0$ and $t_0$ 
such that mild solutions of the two initial value problems satisfy
\begin{equation}
|| u(t,x) - (A_1(t,x)e^{ikx+i\omega t} + A_2(t,x) e^{-ikx-i\omega t}) || 
\leq C\eta = C\varepsilon ^{1/M},
\label{1-12}
\end{equation}
for $t_0 \leq t \leq T_0/\varepsilon $.
\\

If we suppose $A_1 = A_2 : = A$, the system (\ref{1-10b}) is reduced to a single equation
$\partial A/\partial t = \mathcal{Q}A + \varepsilon R(A,A)$; 
the set $\{ A_1 = A_2\}$ is an invariant set of (\ref{1-10b}).
Thus we put $S(A) = R(A,A)$ and consider two initial value problems
\begin{equation}
\displaystyle \frac{\partial u}{\partial t}  = \mathcal{P}u + \varepsilon f(u), \quad
\displaystyle u(0,x) = v_1(\eta x)e^{ikx} + v_1(\eta x)e^{-ikx} 
\label{1-10c}
\end{equation}
and
\begin{equation}
\frac{\partial A}{\partial t} = \mathcal{Q}A + \varepsilon S(A), \quad
\displaystyle A(0,x) = v_1(\eta x).
\label{1-10d} 
\end{equation}
For example if $f(u) = u-u^3$, then $S(A) = A - 3A^3$.
For the equation (\ref{1-10c}), we further suppose that 
\\[0.2cm]
\textbf{(A4)} $P(i\xi) = \overline{P(-i\xi)}$ and $f(\overline{u}) = \overline{f(u)}$.
\\[0.2cm]
That is, $P(x)$ and $f(u)$ are real-valued for $x, u\in \mathbf{R}$.
In particular, if $v_1(\eta x) \in \mathbf{R}$ so that $u(0,x)$ is real-valued,
then a solution $u(t,x)$ is also real-valued.
In the next theorem, $BC^r(\mathbf{R}; \mathbf{R})$ denotes the set of \textit{real}-valued
functions $f$ on $\mathbf{R}$ such that $f(x), f'(x),$ $\cdots ,f^{(r)}(x)$ are bounded uniformly continuous.
\\[0.2cm]
\textbf{Theorem \thedef.}  Suppose $f : BC^r(\mathbf{R}; \mathbf{R}) \to BC^r(\mathbf{R}; \mathbf{R})\, (r\geq 1)$ is $C^2$ such that
the second derivative is locally Lipschitz continuous.
Suppose that there exists a constant $\phi\in \mathbf{R}$ such that $S(\phi) = 0$ and $S'(\phi) < 0$
(that is, $A(t,x) \equiv \phi$ is an asymptotically stable steady state of Eq.(\ref{1-10d})).
If $\varepsilon >0$ is sufficiently small, Eq.(\ref{1-10c}) has a solution of the form
\begin{equation}
u_p(t,x) = \Bigl( \phi + \eta \psi (t, x, \eta) \Bigr) \cdot 2\cos( kx + \omega t).
\end{equation}
The functions $\psi $ and $u_p$ are bounded as $\eta \to 0$ and satisfy
\begin{eqnarray*}
\left\{ \begin{array}{l}
\text{$2\pi / k $-periodic in $x$}\,\, (\text{when}\,\, k\neq 0) \\
\text{$2\pi / \omega $-periodic in $t$\, (when $\omega \neq 0$)}, \\
\text{constant in $t, x$ (when $\omega = 0, k= 0$, respectively)},  \\
\end{array} \right.
\end{eqnarray*}
This $u_p$ is stable in the following sense:
there is a neighborhood $U\subset BC^r(\mathbf{R}; \mathbf{R})$ of $\phi$ in $BC^r(\mathbf{R}; \mathbf{R})$ such that if $v_1\in U$,
then a mild solution $u$ of the initial value problem (\ref{1-10c}) satisfies
$|| u(t, \cdot) - u_p(t, \cdot) || \to 0$ as $t\to \infty$.
\\

The assumption (A4) and a space $BC^r(\mathbf{R}; \mathbf{R})$ means that this theorem holds when every data are real numbers.
The periodic solution $u_p$ is not asymptotically stable toward a complex direction.
These theorems are obtained as special cases of Theorems 4.2 and 4.3 proved in Sec.4.
\\[0.2cm]
\textbf{Example \thedef.} For the Swift-Hohenberg equation, the estimate (\ref{1-4}) immediately follows
from Thm.1.1 by putting $A_2 = \overline{A}_1$ and $v_2 = \overline{v}_1$.
Note that the assumption for the initial value $v_1$ is more relaxed than that given in \cite{Col2}
because we use a mild solution.
To prove the existence of a spatially periodic solution, 
note that the function $S(A)$ is given as $S(A) = A - 3A^3$,
so that $A =\phi = 1/\sqrt{3}$ satisfies the assumptions for Thm.1.2.
Then, it turns out that Eq.(\ref{1-2}) has a stable solution of the form
\begin{equation}
u_p(t,x) = \frac{2}{\sqrt{3}} \cos (kx) + O(\eta),
\end{equation}
which can be obtained directly without using the amplitude equation \cite{Col1}.
\\

The above theorems will be extended to more higher dimensional problems in Sec.4.
\\[0.2cm]
\textbf{Example \thedef.} 
Consider the boundary value problem (\ref{1-14}) with constants $l, d$ and $k$.
Let 
\begin{equation}
L\left(
\begin{array}{@{\,}c@{\,}}
u \\
v
\end{array}
\right) = \left(
\begin{array}{@{\,}c@{\,}}
d(\partial ^2_{x}u + \partial ^2_{y}u) + ku-v \\
\partial ^2_{x}v + \partial ^2_{y}v + u-v
\end{array}
\right)
\label{1-15}
\end{equation}
be the linear operator which defines the unperturbed part.
We suppose that $d$ and $k$ satisfy $(k+d)^2 = 4d$.
Then, the spectrum $\sigma (L)$ is the negative real axis and the origin (See Sec.2),
so that the system (\ref{1-14}) undergoes the Turing instability when $\varepsilon =0$.
The eigenfunction for $0$-eigenvalue satisfying the boundary condition is given by
\begin{equation}
\left(
\begin{array}{@{\,}c@{\,}}
1 \\
(k+d)/2
\end{array}
\right) e^{icy}+\left(
\begin{array}{@{\,}c@{\,}}
1 \\
(k+d)/2
\end{array}
\right) e^{-icy},\quad 
c:= \sqrt{\frac{k-d}{2d}},\, l := \frac{\pi}{c}.
\end{equation}
We will show that the corresponding amplitude equation is given by
\begin{equation}
\frac{\partial A}{\partial t}
 = - \frac{2d^2}{(k+d)(1-d)} \frac{\partial ^4A}{\partial x^4} 
     + \frac{\varepsilon }{1-d}(A - 3A|A|^2), \quad A(0,x) = v_0(\eta x),
\label{1-17}
\end{equation}
where $\eta = \varepsilon ^{1/4}$.
Let us consider a solution of (\ref{1-14}) with the initial condition
\begin{equation}
\left(
\begin{array}{@{\,}c@{\,}}
u(0,x,y) \\
v(0,x,y)
\end{array}
\right) = A(0,x) \left(
\begin{array}{@{\,}c@{\,}}
1 \\
(k+d)/2
\end{array}
\right) e^{icy} + \overline{A(0,x)} \left(
\begin{array}{@{\,}c@{\,}}
1 \\
(k+d)/2
\end{array}
\right) e^{-icy}.
\end{equation}
From theorems shown in Sec.4, it turns out that $u$ is approximately given by
\begin{equation}
u(t,x,y) = A(t,x)e^{icy} + \overline{A(t,x)} e^{-icy} + O(\eta),
\end{equation}
for $t_0\leq t \leq T/\varepsilon $.
Further, the system (\ref{1-14}) proves to have a steady solution of the form
\begin{equation}
\left(
\begin{array}{@{\,}c@{\,}}
u_p(t,x,y) \\
v_p(t,x,y)
\end{array}
\right) = \frac{2}{\sqrt{3}} \cos (cy) \left(
\begin{array}{@{\,}c@{\,}}
1 \\
(k+d)/2
\end{array}
\right) + O(\eta),
\end{equation}
which is periodic in $y$ and constant in $x,t$.
The fourth-order amplitude equation (\ref{1-17}) was also derived by \cite{Ei} in a certain formal way
without error estimates of solutions.
The results in this paper assure that (\ref{1-17}) indeed provides approximate solutions and a steady state.
\\

For both examples, the spectra of the unperturbed linear operators are continuous spectra including the origin.
Thus it is expected that when $\varepsilon $ becomes positive from $0$, 
the spectra get across the imaginary axis and bifurcations occur.
Unfortunately, there are no systematic ways to detect such bifurcations because spectrum on the 
imaginary axis is not discrete.
The reduction proposed in this paper provides a systematic way to detect bifurcations;
a bifurcation problem is reduced to that of the amplitude equation, 
although to investigate the amplitude equation is still difficult in general.

In Sec.2, we will demonstrate that how the amplitude equation is derived by means of the RG method.
The RG method here is one of the singular perturbation methods for differential equations
proposed by Chen, Goldenfeld and Oono \cite{Che}.
In \cite{Chi}, it is proved for ODEs that the RG method unifies classical perturbation methods
such as the multiscaling method, the averaging method, normal forms and so on.
In particular, when the spectrum of an unperturbed linear part is discrete and a center manifold
exists, the RG method is equivalent to the center manifold reduction; the amplitude equation
gives the dynamics on the center manifold.
This paper shows that the RG method and the amplitude equation are still valid even when a spectrum 
is not discrete and a center manifold does not exist.
Even if there are no center manifolds, the amplitude equation provides the dynamics near the 
center subspace and it is useful to study bifurcations of a given PDE.
In particular, Thm.1.2 means that a bifurcation may occur at $\varepsilon =0$;
when $\omega =0$, it is a bifurcation of a steady state and when $\omega \neq 0$,
a $t$-periodic solution appears like as a Hopf bifurcation.

Although results in this paper are partially obtained by many authors for specific problems \cite{Col1,Col2,Sch},
our proof is systematic which is applicable to a wide class of PDEs.
From our proofs, it turns out that reductions of linear operators (reduction of a given differential operator
$\mathcal{P}$ to $\mathcal{Q}$) and that of nonlinearities (reduction of $f(u)$ to $S(u)$) can be done independently.
The reduction of linear operators is described in Prop.3.1 and 3.6, which have the following significant meaning:
the semigroup $e^{\mathcal{P}t}$ generated by $\mathcal{P}$ is approximated by its self-similar part.
Note that the evolution equation $\dot{u} = \mathcal{Q}u$ has the self-similar structure in the sense that it is 
invariant under the transformation
\begin{equation}
(t, x) \mapsto (c^M t, cx), \quad c\in \mathbf{R},
\end{equation}
see (\ref{1-8}).
Then, Prop.3.1 implies that a \textit{non}-self-similar part of $e^{\mathcal{P}t}$ decays to zero as $t\to \infty$.
In other words, if we apply the above transformation repeatedly, then a non-self-similar part decays to zero,
while a self-similar part  survives because it is invariant under the transformation.
This self-similar part defines the linear operator $\mathcal{Q}$.
Such a technique to obtain a self-similar structure is also known as the renormalization group method 
in statistical mechanics.

%%%%%%%%%%%%%%%%%%%%%%%%%%%%%%%%%%%%%%%%%%%%%%%%%%%%%%%%%%%%%%%%%%%%%%%%%%%%%%%%%%%%%%%%%%%%%%%%%%
%%%%%%%%%%%%%%%%%%%%%%%%%%%%%%%%%%%%%%%%%%%%%%%%%%%%%%%%%%%%%%%%%%%%%%%%%%%%%%%%%%%%%%%%%%%%%%%%%%

\section{The renormalization group method}

In this section, we demonstrate how the amplitude equation for (\ref{1-1}) is obtained
by the RG method with examples.
The RG method is a formal way to find amplitude equations and the results in this section will not
be used in later sections.
Although we only consider parabolic-type PDEs, the RG method is applicable to a more large class of PDEs
and it has advantages over the multiscaling method \cite{Che}.
See \cite{Chi} for the RG method for ODEs.

Let us consider the Swift-Hohenberg equation (\ref{1-2}).
We expand a solution as $u = u_0 + \varepsilon u_1 + O(\varepsilon ^2)$.
The zero-th order term $u_0$ satisfies the linear equation
\begin{equation}
\frac{\partial u_0}{\partial t} = -(\partial _x^2 + k^2)^2u_0.
\end{equation}
We are interested in the dynamics near the center subspace.
The spectrum of $-(\partial _x^2 + k^2)^2$ intersects with the imaginary axis at the origin,
and the corresponding eigenfunctions are $e^{ikx}$ and $e^{-ikx}$.
Thus we consider the solution of the form
\begin{equation}
u_0(x) = Ae^{ikx} + Be^{-ikx},
\end{equation}
where $A, B\in \mathbf{C}$ are constants.
Then, the first order term $u_1$ satisfies the inhomogeneous linear equation
\begin{eqnarray}
& & \frac{\partial u_1}{\partial t} 
= -(\partial _x^2 + k^2)^2u_1 + u_0-u_0^3 \nonumber \\
&=& -(\partial _x^2 + k^2)^2u_1
 + Ae^{ikx} + Be^{-ikx}-(A^3e^{3ikx} + 3A^2Be^{ikx}+3AB^2 e^{-ikx}+B^3e^{-3ikx}). \quad
\label{2-3}
\end{eqnarray}
Since the factors $e^{\pm ikx}$ in the inhomogeneous terms are eigenfunctions of the operator
$-(\partial _x^2 + k^2)^2$, it is expected that a solution of this equation includes secular terms
which diverge in $t$ and $x$.
To find secular terms arising from the factor $e^{ikx}$, we consider the equation
\begin{eqnarray}
\frac{\partial u_1}{\partial t} 
= -(\partial _x^2 + k^2)^2u_1 + (A-3A^2B)e^{ikx}
\label{2-4}
\end{eqnarray}
instead of Eq.(\ref{2-3}).
We assume a special solution of the form
\begin{eqnarray*}
u_1 = (\mu_1 t^{\beta_1} + \mu_2 x^{\beta_2}) e^{ikx}.
\end{eqnarray*}
Substituting this into Eq.(\ref{2-4}), we obtain $\beta_1 =1, \, \beta_2 = 2$,
and $\mu_1, \mu_2$ prove to satisfy the relation
\begin{equation}
\mu_1 = 8k^2 \mu_2 + A-3A^2B.
\label{2-5}
\end{equation} 
Secular terms corresponding to the factor $e^{-ikx}$ are obtained in the same way.
Therefore, a special solution of Eq.(\ref{2-3}) including secular terms are given by
\begin{eqnarray*}
u_1 = (\mu_1 t + \mu_2 x^2) e^{ikx} + (\tilde{\mu}_1 t + \tilde{\mu}_2 x^2) e^{-ikx}
 - \frac{A^3}{64}e^{3ikx}  - \frac{B^3}{64}e^{-3ikx},
\end{eqnarray*}
where $\tilde{\mu}_1$ and $\tilde{\mu}_2$ satisfy the same relation as (\ref{2-5}).
Thus we obtain
\begin{eqnarray*}
u &=& Ae^{ikx} + Be^{-ikx} \\
& + & \varepsilon \left( (\mu_1 t + \mu_2 x^2) e^{ikx} + (\tilde{\mu}_1 t + \tilde{\mu}_2 x^2) e^{-ikx}
 + (\mathrm{nonsecular}) \right) + O(\varepsilon ^2).
\end{eqnarray*}
In what follows, we omit to write down nonsecular terms and $O(\varepsilon ^2)$-terms 
which will not be used later.
To remove the secular terms, we introduce dummy parameters $\tau$ and $X$,
and rewrite the above $u$ as
\begin{eqnarray*}
u &=& Ae^{ikx} + Be^{-ikx} + (\varepsilon \mu_1 \tau + \varepsilon \mu_2 X^2)e^{ikx} 
          + (\varepsilon \tilde{\mu}_1 \tau + \varepsilon \tilde{\mu}_2 X^2)e^{-ikx} \\
& + & \varepsilon \left( \mu_1 (t-\tau) + \mu_2 (x^2-X^2) \right) e^{ikx} 
 +    \varepsilon \left( \tilde{\mu}_1 (t-\tau) + \tilde{\mu}_2 (x^2-X^2) \right) e^{-ikx}.
\end{eqnarray*}
Now terms $\varepsilon \mu_1 \tau + \varepsilon \mu_2 X^2$ and 
$\varepsilon \tilde{\mu}_1 \tau + \varepsilon \tilde{\mu}_2 X^2$ are \textit{renormalized} into 
the constants $A$ and $B$, respectively.
Thus we rewrite $u$ as
\begin{eqnarray*}
u &=& A(\tau, X)e^{ikx} + B(\tau, X)e^{-ikx} \\
& + & \varepsilon \left( \mu_1 (t-\tau) + \mu_2 (x^2-X^2) \right) e^{ikx} 
 +    \varepsilon \left( \tilde{\mu}_1 (t-\tau) + \tilde{\mu}_2 (x^2-X^2) \right) e^{-ikx}.
\end{eqnarray*}
Putting $\tau = t$ and $X = x$ provides
\begin{eqnarray*}
u(t,x) = A(t, x)e^{ikx} + B(t, x)e^{-ikx},
\end{eqnarray*}
which seems to give an approximate solution if $A(t, x)$ and $B(t, x)$ are appropriately defined.
Since $u$ is independent of dummy parameters $\tau $ and $X$, we require that the equation
\begin{eqnarray*}
\frac{\partial u}{\partial \tau} = \frac{\partial ^2u}{\partial X^2} = 0
\end{eqnarray*}
holds, which is called the \textit{RG equation}.
This yields
\begin{equation}
\left\{ \begin{array}{l}
\displaystyle \frac{\partial u}{\partial \tau}\Bigl|_{\tau = t, X=x}
 = \left( \frac{\partial A}{\partial t} - \varepsilon \mu_1  \right)e^{ikx}
   + \left( \frac{\partial B}{\partial t} - \varepsilon \tilde{\mu}_1  \right)e^{-ikx} = 0, \\[0.3cm]
\displaystyle \frac{\partial^2 u}{\partial X^2}\Bigl|_{\tau = t, X=x}
 = \left( \frac{\partial^2 A}{\partial x^2} - 2\varepsilon \mu_2  \right)e^{ikx}
   + \left( \frac{\partial^2 B}{\partial x^2} - 2\varepsilon \tilde{\mu}_2  \right)e^{-ikx} = 0.
\end{array} \right.
\end{equation}
Since $\mu_1$ and $\mu_2$ satisfy (\ref{2-5}), we obtain
\begin{equation}
\frac{\partial A}{\partial t} = \varepsilon \mu_1 = 
   8k^2 \varepsilon \mu_2 + \varepsilon (A-3A^2B) 
    = 4k^2 \frac{\partial ^2 A}{\partial x^2} + \varepsilon (A-3A^2B) .
\end{equation}
Similarly, $B$ satisfies $\partial B/\partial t = 4k^2 \partial_x^2 B + \varepsilon (B-3AB^2)$.
If we suppose that $A = \overline{B}$ to obtain a real-valued solution,
the Ginzburg-Landau equation (\ref{1-3}) is obtained.
\\

Next, let us derive the amplitude equation of the system (\ref{1-14}).
The dispersion relation of the unperturbed operator $L$ (\ref{1-15}) is
\begin{eqnarray}
& & \det \left(
\begin{array}{@{\,}cc@{\,}}
\lambda +d(\xi_1^2 + \xi_2^2) - k & 1 \\
-1 & \lambda +\xi_1^2 + \xi_2^2 + 1
\end{array}
\right) \nonumber \\ 
&=& 
\lambda ^2 + (d\xi^2 + \xi^2 + 1-k) \lambda + (d\xi^2 - k)(\xi^2 + 1) +1 = 0,
\label{2-8}
\end{eqnarray}
where we put $\xi^2 = \xi_1^2 + \xi_2^2$.
Let $\lambda _{\pm}(\xi)$ be two roots of (\ref{2-8}).
Then, the spectrum of $L$ is given by $\sigma (L) = \lambda _+(\mathbf{R}) \cup \lambda _-(\mathbf{R})$.
Suppose that Eq.(\ref{1-14}) undergoes the Turing instability at $\varepsilon =0$,
so that $\sigma (L) = \mathbf{R}_{\leq 0}$.
It is easy to verify that this is true if and only if 
\begin{eqnarray*}
0<d<k<1,\,\, (k+d)^2 = 4d.
\end{eqnarray*}
In particular, one of $\lambda _{\pm}(\xi)$ satisfies $\lambda _{\pm}(c) = 0$, where $c^2 = (k-d)/2d$.
For any $(\xi_1, \xi_2)$ satisfying $ \xi_1^2 + \xi_2^2 = c^2$,
\begin{eqnarray*}
\left(
\begin{array}{@{\,}c@{\,}}
1 \\
(k+d)/2
\end{array}
\right) e^{i\xi_1 x + i\xi_2 y}
\end{eqnarray*}
is an eigenfunction of $L$ associated with $\lambda =0$.
Because of the boundary condition in (\ref{1-14}), we choose $e^{icy}$ and $e^{-icy}$.
Thus we expand a solution of (\ref{1-14}) as
\begin{equation}
\left(
\begin{array}{@{\,}c@{\,}}
u \\
v
\end{array}
\right) = A \left(
\begin{array}{@{\,}c@{\,}}
1 \\
(k+d)/2
\end{array}
\right) e^{icy} + B\left(
\begin{array}{@{\,}c@{\,}}
1 \\
(k+d)/2
\end{array}
\right) e^{-icy} + \varepsilon \left(
\begin{array}{@{\,}c@{\,}}
u_1 \\
v_1
\end{array}
\right) + O(\varepsilon ^2).
\end{equation}
Put $A = \overline{B}$ for simplicity.
Then, $(u_1, v_1)$ satisfies the equation
\begin{equation}
\frac{\partial }{\partial t}\left(
\begin{array}{@{\,}c@{\,}}
u_1 \\
v_1
\end{array}
\right) = L \left(
\begin{array}{@{\,}c@{\,}}
u_1 \\
v_1
\end{array}
\right) + \left(
\begin{array}{@{\,}c@{\,}}
A e^{icy} - A^3 e^{3icy} - 3A|A|^2e^{icy} + c.c. \\
0
\end{array}
\right),
\label{2-10}
\end{equation}
where $c.c.$ denotes the complex conjugate.
We find secular terms of the form
\begin{equation}
u_1 = (\mu_1 t + \mu_2 x^2 + \mu_3 x^4) e^{icy},
\quad v_1 = \frac{k+d}{2}(\tilde{\mu}_1 t + \tilde{\mu}_2 x^2 + \tilde{\mu}_3 x^4) e^{icy}.
\end{equation}
Substituting them into (\ref{2-10}), we obtain
\begin{equation}
\left\{ \begin{array}{l}
\tilde{\mu}_1 = \mu_1, \\
\tilde{\mu}_3 = \mu_3,  \\
\tilde{\mu}_2 = \frac{1}{2}\mu_1,
\end{array} \right.
\quad \left\{ \begin{array}{l}
\mu_1 = 2d\mu_2 + A - 3A|A|^2, \\[0.2cm]
\displaystyle 0 = 12d \mu_3 + k\mu_2 - \frac{k+d}{2}\tilde{\mu}_2 - dc^2 \mu_2.  \\
\end{array} \right.
\label{2-11}
\end{equation}
Then, a formal solution is given as
\begin{eqnarray*}
u = Ae^{icy} + \varepsilon (\mu_1 t + \mu_2 x^2 + \mu_3 x^4) e^{icy}
 + c.c. + (\mathrm{nonsecular}) + O(\varepsilon ^2).
\end{eqnarray*}
Introducing dummy parameters $\tau, X$ and renormalizing, we rewrite this equation as
\begin{eqnarray*}
u \!=\! A(\tau, X)e^{icy} \!+\! \varepsilon \left( \mu_1 (t\!-\!\tau) + 
         \mu_2 (x^2\!-\!X^2) + \mu_3 (x^4\!-\!X^4) \right) e^{icy}
 \!+ \!c.c.\! + \!(\mathrm{nonsecular})\! + \! O(\varepsilon ^2).
\end{eqnarray*}
Since $u$ is independent of $\tau $ and $X$, we require
\begin{equation}
\left\{ \begin{array}{l}
\displaystyle \frac{\partial u}{\partial \tau}\Bigl|_{\tau = t, X=x}
 = \left( \frac{\partial A}{\partial t} - \varepsilon \mu_1 \right) e^{icy} + c.c. = 0,  \\
\displaystyle \frac{\partial^4 u}{\partial X^4}\Bigl|_{\tau = t, X=x}
 = \left( \frac{\partial^4 A}{\partial X^4} - 24 \varepsilon \mu_3 \right) e^{icy} + c.c. = 0.  \\
\end{array} \right.
\label{2-12}
\end{equation}
Finally, Eqs.(\ref{2-11}) and (\ref{2-12}) provide the amplitude equation (\ref{1-17})
by eliminating $\mu_1, \mu_2$ and $\mu_3$.

%%%%%%%%%%%%%%%%%%%%%%%%%%%%%%%%%%%%%%%%%%%%%%%%%%%%%%%%%%%%%%%%%%%%%%%%%%%%%%%%%%%%%%%%%%%%%%%%%%
%%%%%%%%%%%%%%%%%%%%%%%%%%%%%%%%%%%%%%%%%%%%%%%%%%%%%%%%%%%%%%%%%%%%%%%%%%%%%%%%%%%%%%%%%%%%%%%%%%

\section{Reduction of a linear semigroup}

For Eq.(\ref{1-1}), reductions of the linear unperturbed part $\mathcal{P}u$
and the perturbation term $f(u)$ can be done independently.
In this section, we give a reduction of the linear part.

%%%%%%%%%%%%%%%%%%%%%%%%%%%%%%%%%%%%%%%%%%%%%%%%%%%%%%%%%%%%%%%%%%%%%%%%%%%%%%%%%%%%%%%%%%%%%%%%%%

\subsection{One dimensional case}

We start with the simple case $u\in \mathbf{C}$ and $x\in \mathbf{R}^d$.
Put $x = (x_1, \cdots ,x_d)$ and $\alpha = (\alpha _1, \cdots ,\alpha _d)$,
where $\alpha$ denotes a multi-index as usual:
$x^\alpha  = (x_1^{\alpha _1} , \cdots ,x_d^{\alpha _d})$ and $|\alpha | = \alpha _1+ \cdots +\alpha _d$.
Let $P(x) = \sum^q_{|\alpha |=0} a_\alpha x^\alpha $ be a polynomial of degree $q$ and 
$\mathcal{P}:= P(\partial _1 , \cdots ,\partial _d)$ a differential operator on $\mathbf{R}^d$,
where $\partial _j$ denotes the derivative with respect to $x_j$.
We make the following assumptions.
\\[0.2cm]
\textbf{(B1)} $\mathrm{Re}[P(i\xi)] \leq 0$ for any $\xi \in \mathbf{R}^d$.
\\
\textbf{(B2)} There exist $\omega \in \mathbf{R},\, k\in \mathbf{R}^d$ and an integer $M$ such that
\begin{eqnarray*}
& & P(ik) = i\omega, \\
& & \frac{\partial ^\alpha P}{\partial x^\alpha }(ik) = 0, \,\,\text{for any $\alpha $ such that}
\,\, |\alpha | = 1, \cdots  ,M-1, \\
& & \frac{\partial ^\alpha P}{\partial x^\alpha }(ik) \neq 0,\,\, \text{for some $\alpha $ such that}
\,\, |\alpha | = M.
\end{eqnarray*}
\textbf{(B3)} Define $Q(x)$ and $\mathcal{Q}$ by
\begin{equation}
Q(x) = \sum_{|\alpha | = M} \frac{1}{(\alpha _1!)\cdots (\alpha _d !)} 
\frac{\partial ^\alpha P}{\partial x^\alpha }(ik) x^\alpha , \quad 
\mathcal{Q} = Q(\partial _1 , \cdots ,\partial _d).
\end{equation}
Then, both of $\mathcal{P}$ and $\mathcal{Q}$ are elliptic in the sense that
there exist $c_1,c_2 > 0$ such that $\mathrm{Re}[P(i\xi )] < -c_2|\xi|^2$ and 
$\mathrm{Re}[Q(i\xi)] < -c_2|\xi|^2$ hold for $|\xi|\geq c_1$.
\\

Put $B^r = BC^r(\mathbf{R}^d; \mathbf{C})$, a Banach space
of complex-valued bounded uniformly continuous functions on $\mathbf{R}^d$ up to the $r$-th derivative.
In the next propositions, $|| \cdot || = || \cdot ||_r$ denotes the standard supremum norm on $B^r$.
Consider two initial value problems:
\begin{subnumcases}
{}
\displaystyle \frac{\partial u}{\partial t}  = \mathcal{P}u, \quad
\displaystyle u(0,x) = v_0(x)e^{ikx}, \label{3-2a}\\[0.2cm]
\displaystyle \frac{\partial A}{\partial t} = \mathcal{Q}A. \quad
\displaystyle A(0,x) = v_0(x),\label{3-2b}
\end{subnumcases}
where $kx = k_1 x_1 + \cdots + k_dx_d$, and a similar notation will be used in the sequel.
Because of (B3), $\mathcal{P}$ and $\mathcal{Q}$ 
generate $C^0$-semigroups $e^{\mathcal{P}t}$ and $e^{\mathcal{Q}t}$ on $B^r$, respectively.
Thus solutions of the above problems are written as $e^{\mathcal{P}t}(e^{ikx} v_0)$
and $e^{\mathcal{Q}t}v_0$, respectively.
\\[0.2cm]
\textbf{Proposition \thedef.}
Suppose (B1) to (B3) and $r\geq 0$.
There exists a constant $C_1>0$ such that the inequality
\begin{equation}
|| e^{\mathcal{P}t}(e^{ikx} v_0) - e^{i\omega t + ikx}e^{\mathcal{Q}t}v_0 ||_r \leq C_1 t^{-1/M} || v_0 ||_r
\end{equation}
holds for any $t>0$ and $v_0 \in BC^r(\mathbf{R}^d; \mathbf{C}) $.

For the main theorems in this paper, we need the following perturbative problem
\begin{subnumcases}
{}
\displaystyle \frac{\partial u}{\partial t}  = \mathcal{P}u , \quad
\displaystyle u(0,x) = v_0(\eta x)e^{ikx}, \label{3-4a}\\[0.2cm]
\displaystyle \frac{\partial A}{\partial t} = \mathcal{Q}A, \quad
\displaystyle A(0,x) = v_0(\eta x),\label{3-4b}
\end{subnumcases}
where $\eta = \varepsilon ^{1/M}$ and $\varepsilon >0$ is a small parameter.
\\[0.2cm]
\textbf{Proposition \thedef.}
Suppose (B1) to (B3) and $r\geq 1$.
For any $\varepsilon >0$ and $t_0 > 0$, there exists a positive number $C_1 = C_1(t_0)$ such that the inequality
\begin{equation}
|| e^{\mathcal{P}t}(e^{ikx} \hat{v}_0) - e^{i\omega t + ikx}e^{\mathcal{Q}t}\hat{v}_0 || _r
\leq \eta C_1 || v_0 ||_r
\label{3-5}
\end{equation}
holds for $t \geq t_0$ and $v_0 \in BC^r(\mathbf{R}^d; \mathbf{C})$, where $\hat{v}_0(\cdot) := v_0(\eta \,\cdot)$.
\\[0.2cm]
\textbf{Proof of Prop.3.1.} 
By putting $u = e^{i\omega t}w$, Eq.(\ref{3-2a}) is rewritten as 
$\partial w/\partial t = (\mathcal{P} - i\omega )w$.
Then, the operator $\mathcal{P}-i\omega $ satisfies (B1) to (B3) with $\omega =0$.
Hence, it is sufficient to prove the proposition for $\omega =0$.

Two solutions are given by
\begin{eqnarray*}
A(t,x) = \frac{1}{(2\pi)^d} \int\! v_0(y+x) \int\! e^{-iy\xi} e^{Q(i\xi)t} d\xi dy  
\end{eqnarray*}
and
\begin{eqnarray}
u(t,x) &=& \frac{1}{(2\pi)^d} \int\! v_0(y+x)e^{ik(y+x)} \int\! e^{-iy\xi} e^{P(i\xi)t} d\xi dy \nonumber \\
&=& \frac{e^{ikx}}{(2\pi)^d} \int\! v_0(y+x) \int\! e^{-iy\xi} e^{P(i\xi + ik)t} d\xi dy,
\label{3-5b}
\end{eqnarray}
respectively.
Thus we obtain
\begin{eqnarray*}
u(t,x) - e^{ikx}A(t,x) 
= \frac{e^{ikx}}{(2\pi)^d} \int\! v_0(y+x) \int\! e^{-iy\xi} e^{Q(i\xi)t}
\left( e^{P(i\xi + ik)t - Q(i\xi)t}  - 1\right) d\xi dy.
\end{eqnarray*}
Put $\tau = t^{-1/M}$.
Changing variables $\xi \mapsto \tau \xi,\, y\mapsto y/\tau$ yields
\begin{eqnarray*}
u(t,x) - e^{ikx}A(t,x) 
= \frac{e^{ikx}}{(2\pi)^d} \int\! v_0(y/\tau +x) \int\! e^{-iy\xi} e^{Q(i\xi)}
\left( e^{P(i\tau \xi + ik)/\tau^M - Q(i\xi)}  - 1\right) d\xi dy.
\end{eqnarray*}
Due to the assumption (B2), we have
\begin{eqnarray*}
g(\xi, \tau) &:=& P(i\tau \xi + ik)/\tau^M - Q(i\xi) \\
&=& \frac{1}{\tau^M} \sum^q_{|\alpha |=0} \frac{1}{(\alpha _1!)\cdots (\alpha _d !)}
\frac{\partial ^\alpha P}{\partial x^\alpha }(ik) i^{|\alpha |}\tau^{|\alpha |} \xi^\alpha 
 - \sum_{|\alpha | = M}\frac{1}{(\alpha _1!)\cdots (\alpha _d !)}
\frac{\partial ^\alpha P}{\partial x^\alpha }(ik) i^{M}\xi^\alpha \\
&=& \sum^q_{|\alpha | = M+1}\frac{1}{(\alpha _1!)\cdots (\alpha _d !)}
\frac{\partial ^\alpha P}{\partial x^\alpha }(ik) i^{|\alpha |}\xi^\alpha \cdot \tau^{|\alpha |-M}. 
\end{eqnarray*}
Note that $g\sim O(\tau)$ as $\tau \to 0$.
In particular, there exists $0<\theta <1$ such that 
\begin{eqnarray*}
e^{g(\xi, \tau)} - 1 = \tau \frac{\partial g}{\partial \tau}(\xi, \theta \tau)e^{g(\xi, \theta \tau)}.
\end{eqnarray*}
This provides
\begin{equation}
\left\{ \begin{array}{l}
\displaystyle u(t,x) - e^{ikx}A(t,x) = \tau \frac{e^{ikx}}{(2\pi)^d} \int\! v_0(y/\tau +x)G(y, \tau) dy, \\[0.3cm]
\displaystyle G(y, \tau) := \int\! e^{-iy\xi} e^{Q(i\xi)} 
\frac{\partial g}{\partial \tau}(\xi, \theta \tau)e^{g(\xi, \theta \tau)} d\xi.   \\
\end{array} \right.
\label{3-6}
\end{equation}
Because of (B3), $G(y, \tau)$ exists for each $\tau \geq 0$ and $y\in \mathbf{R}$.
Since $g$ is polynomial in $\tau$,
there exist $\tau_0$ and $D_1 = D_1(\tau_0)$ such that $|G(y, \tau)| \leq D_1$ holds 
for $0\leq \tau \leq \tau_0$ and $y\in [-1,1]^d$.
Next, since the integrand in the definition of $G(y, \tau)$ is smooth in $\xi$,
$G(y, \tau)$ is rapidly decreasing in $y$ due to the property of the Fourier transform.
Indeed, by using integration by parts, it is easy to verify that there exists $D_2 = D_2(\tau_0)$
such that $|G(y, \tau)| \leq D_2 (y_1 \cdots  y_d)^{-2}$ holds for $0\leq \tau \leq \tau_0$ and $y\notin [-1,1]^d$.
This provides
\begin{eqnarray*}
|u(t,x) - e^{ikx}A(t,x)| 
&\leq & \frac{\tau}{(2\pi)^d} \int\! |v_0(y/\tau + x)| \cdot |G(y,\tau)| dy \\
&\leq &  \frac{\tau}{(2\pi)^d} D_1 || v_0 || + 
         \frac{\tau}{(2\pi)^d} \int_{y\notin [-1,1]^d}\! \frac{D_2}{y_1^2 \cdots  y_d^2} dy \cdot|| v_0 ||.
\end{eqnarray*}
This proves that 
\begin{equation}
\sup_{x\in \mathbf{R}^d} |u(t,x) - e^{ikx}A(t,x)|\leq \tau D_3 || v_0 ||
\end{equation}
for some $D_3 > 0$ when $0\leq \tau \leq \tau_0$.
To estimate the derivatives, note that Eq.(\ref{3-5b}) is rewritten as
\begin{eqnarray*}
e^{-ikx}u(t,x) = \frac{1}{(2\pi)^d} \int\! v_0(y) \int\! e^{-i(y-x)\xi} e^{P(i\xi + ik)t} d\xi dy,
\end{eqnarray*}
and similarly for $A(t,x)$.
Hence, the derivative is given as
\begin{eqnarray*}
\left\{ \begin{array}{l}
\displaystyle \frac{\partial ^\alpha }{\partial x^\alpha } \left(e^{-ikx}u(t,x) - A(t,x)  \right)
       = \tau \frac{1}{(2\pi)^d} \int\! v_0(y/\tau +x)G_\alpha (y, \tau) dy, \\[0.3cm]
\displaystyle G_\alpha (y, \tau) := \int\! (i\tau \xi)^\alpha e^{-iy\xi} e^{Q(i\xi)} 
\frac{\partial g}{\partial \tau}(\xi, \theta \tau)e^{g(\xi, \theta \tau)} d\xi.   \\
\end{array} \right.
\end{eqnarray*}
By the same way as above, we can show that this derivative is of $O(\tau)$ uniformly in $x$.
Hence, the inequality
\begin{equation}
|| u(t,x) - e^{ikx}A(t,x) || \leq \tau D_3 || v_0 || = t^{-1/M} D_3 || v_0 ||
\label{3-9}
\end{equation}
holds with respect to the norm of $B^r$ for some $D_3>0$ and any $t \geq \tau_0^{-M}$.
On the other hand, since $\mathcal{P}$ and $\mathcal{Q}$ generate $C^0$-semigroups on $B^r$,
there exists $D_4 > 0$ such that $|| u(t,x) - e^{ikx}A(t,x) || \leq D_4|| v_0 ||$ for 
$0\leq t\leq \tau_0^{-M}$.
This and Eq.(\ref{3-9}) prove Prop.3.1 (for $\omega  = 0$). \hfill $\blacksquare$
\\[0.2cm]
\textbf{Proof of Prop.3.2.}
In this case, solutions satisfy
\begin{eqnarray*}
u(t,x) - e^{ikx}A(t,x) 
= \tau \frac{e^{ikx}}{(2\pi)^d} \int\! v_0(\eta y/\tau + \eta x)G(y, \tau) dy,
\end{eqnarray*}
where $\tau = t^{-1/M}$ and $G$ is defined by (\ref{3-6}) as before.
Since $v_0\in B^{r}\, (r\geq 1)$, there exists $0< \theta _1 < 1$ such that it is expanded as
\begin{eqnarray*}
u(t,x) - e^{ikx}A(t,x) 
&=& \tau \frac{e^{ikx}}{(2\pi)^d} \int\! \left( v_0(\eta x) + \sum^d_{j=1} \frac{\partial v_0}{\partial x_j} 
         (\eta x + \theta _1 \eta y/\tau ) \frac{\eta}{\tau} y_j \right)G(y, \tau) dy \\
&=& \eta \frac{e^{ikx}}{(2\pi)^d} \int\! \sum^d_{j=1} \frac{\partial v_0}{\partial x_j} 
        (\eta x + \theta _1 \eta y/\tau ) y_j G(y, \tau) dy,
\end{eqnarray*}
where we used the fact $\int\! G(y, \tau) dy = 0$.
The rest of the proof is the same as that of Prop.3.1. 
\hfill $\blacksquare$
\\

If the polynomial $P(x)$ has no symmetries, the assumptions (B2),(B3) seem to be strong;
for example, if $d=2$ and 
\begin{eqnarray*}
\frac{\partial ^2 P}{\partial x_1^2}(ik) \neq 0, \quad 
\frac{\partial ^2 P}{\partial x_1x_2}(ik) = \frac{\partial ^2 P}{\partial x_2^2}(ik) = 0,
\end{eqnarray*}
then $\mathcal{Q}$ is not elliptic.
To relax the assumptions, fix an integer $D$ such that $1\leq D \leq d$.
We denote $x \in \mathbf{R}^d$ as $x = (\hat{x}_1, \hat{x}_2)$ with $\hat{x}_1 = (x_1, \cdots , x_D)$
and $\hat{x}_2 = (x_{D+1} , \cdots ,x_d)$.
Accordingly, a multi-index $\alpha $ is also denoted as $\alpha  = (\beta, \gamma )$.
Instead of (B2) and (B3), we suppose that 
\\[0.2cm]
\textbf{(B2)$_D$}\, there exist $\omega \in \mathbf{R},\, k\in \mathbf{R}^d$ and an integer $M$ such that
\begin{eqnarray*}
& & P(ik) = i\omega, \\
& & \frac{\partial ^\beta P}{\partial \hat{x}_1^\beta }(ik) = 0, \,\,\text{for any $\beta $ such that}
\,\, |\beta | = 1, \cdots  ,M-1, \\
& & \frac{\partial ^\beta P}{\partial \hat{x}_1^\beta }(ik) \neq 0,\,\, \text{for some $\beta $ such that}
\,\, |\beta | = M.
\end{eqnarray*}
\textbf{(B3)$_D$} Define $Q(x)$ and $\mathcal{Q}$ by
\begin{equation}
Q(x) = Q(\hat{x}_1, 0) = \sum_{|\beta | = M} \frac{1}{(\beta _1!)\cdots (\beta _D !)} 
\frac{\partial ^\beta P}{\partial \hat{x}_1^\beta }(ik) \hat{x}_1^\beta , \quad 
\mathcal{Q} = Q(\partial _1 , \cdots ,\partial _D, 0 ,\cdots ,0).
\end{equation}
Then, both of $\mathcal{P}$ and $\mathcal{Q}$ are elliptic in the sense that
there exist $c_1,c_2 > 0$ such that $\mathrm{Re}[P(i\xi )] \!<\! -c_2|\xi|^2$ and 
$\mathrm{Re}[Q(i\hat{\xi}_1)] \!<\! -c_2|\hat{\xi}_1|^2$ hold for $|\xi|, |\hat{\xi}_1| \geq c_1$,
where $\hat{\xi}_1 = (\xi_1, \cdots ,\xi_D)$.
\\

When $D = d$, these assumptions are reduced to (B2) and (B3) before.
The assumption (B3)$_D$ implies that $\mathcal{Q}$ is an elliptic operator on $\mathbf{R}^D$
although it is not on $\mathbf{R}^d$.
Consider two initial value problems:
\begin{subnumcases}
{}
\displaystyle \frac{\partial u}{\partial t}  = \mathcal{P}u, \quad
\displaystyle u(0,x) = v_0(\hat{x}_1)e^{ikx}, \label{3-11a}\\[0.2cm]
\displaystyle \frac{\partial A}{\partial t} = \mathcal{Q}A. \quad
\displaystyle A(0,x) = v_0(\hat{x}_1). \label{3-11b}
\end{subnumcases}
Note that $v_0$ depends only on $\hat{x}_1 = (x_1, \cdots  ,x_D)$.
In particular, Eq.(\ref{3-11b}) can be regarded as a parabolic equation on $BC^r(\mathbf{R}^{D}, \mathbf{C})$,
while Eq.(\ref{3-11a}) is a parabolic equation on $BC^r(\mathbf{R}^{d}, \mathbf{C})$.
We also consider the perturbative problem
\begin{subnumcases}
{}
\displaystyle \frac{\partial u}{\partial t}  = \mathcal{P}u , \quad
\displaystyle u(0,x) = v_0(\eta \hat{x}_1)e^{ikx}, \label{3-12a}\\[0.2cm]
\displaystyle \frac{\partial A}{\partial t} = \mathcal{Q}A, \quad
\displaystyle A(0,x) = v_0(\eta \hat{x}_1),\label{3-12b}
\end{subnumcases}
where $\eta = \varepsilon ^{1/M}$ and $\varepsilon >0$ is a small parameter.
\\[0.2cm]
\textbf{Proposition \thedef.}
Under the assumptions (B1), (B2)$_D$ and (B3)$_D$, Prop.3.1 holds for Eqs.(\ref{3-11a}),(\ref{3-11b}),
and Prop.3.2 holds for Eqs.(\ref{3-12a}),(\ref{3-12b}).
\\[0.2cm]
\textbf{Proof.}
For Eq.(\ref{3-11a}), $u(t,x)$ is given as
\begin{eqnarray*}
u(t,x) &=& \frac{e^{ikx}}{(2\pi)^d} \int\! v_0(\hat{y}_1+\hat{x}_1) \int\! e^{-iy\xi} e^{P(i\xi + ik)t} d\xi dy \\
&=& \frac{e^{ikx}}{(2\pi)^d} \int\! v_0(\hat{y}_1+\hat{x}_1) 
    \int\! e^{-i\hat{y}_1 \hat{\xi}_1 - i \hat{y}_2 \hat{\xi}_2} 
    e^{P(i\hat{\xi}_1 + i\hat{k}_1 , i\hat{\xi}_2 + i\hat{k}_2 )t} d\hat{\xi}_1d\hat{\xi}_2d\hat{y}_1d\hat{y}_2.
\end{eqnarray*}
To calculate this, we need the next lemma.
\\[0.2cm]
\textbf{Lemma \thedef.} Let $\mathcal{S}$ be a space of $C^\infty$ rapidly decreasing functions on 
$\mathbf{R}^{d-D}$ (Schwartz space). For any $f\in \mathcal{S}$, we have
\begin{equation}
\int\!\! \int e^{-i \hat{y}_2 \hat{\xi}_2} f(\hat{\xi}_2) d\hat{\xi}_2d\hat{y}_2 =
(2\pi)^{d-D} f(0).  
\end{equation}
\textbf{Proof.} Let $\mathcal{S}'$ be a dual space of $\mathcal{S}$.
For the pairing of $\mathcal{S}'$ and $\mathcal{S}$, we use a bracket $\langle \,\, \,,\, \,\, \rangle$.
Let $\mathcal{F}$ be the Fourier transform.
Then,
\begin{eqnarray*}
\int\!\! \int e^{-i \hat{y}_2 \hat{\xi}_2} f(\hat{\xi}_2) d\hat{\xi}_2d\hat{y}_2 
&=& (2\pi)^{(d-D)/2} \int \! \mathcal{F}[f](\hat{y}_2) d\hat{y}_2 \\
&=& (2\pi)^{(d-D)/2} \langle 1 \,,\, \mathcal{F}[f] \rangle \\
&=& (2\pi)^{(d-D)/2} \langle \mathcal{F}[1] \,,\, f \rangle \\
&=& (2\pi)^{d-D} \langle \delta \,,\, f \rangle = (2\pi)^{d-D}f(0),
\end{eqnarray*}
where $\delta $ is the Dirac delta. \hfill $\blacksquare$
\\

Due to this lemma, we obtain
\begin{eqnarray*}
u(t,x) = \frac{e^{ikx}}{(2\pi)^D} \int\! v_0(\hat{y}_1+\hat{x}_1) 
    \int\! e^{-i\hat{y}_1 \hat{\xi}_1} e^{P(i\hat{\xi}_1 + i\hat{k}_1 , i\hat{k}_2 )t} d\hat{\xi}_1d\hat{y}_1.
\end{eqnarray*}
Since Eq.(\ref{3-11b}) is an equation on $BC^r(\mathbf{R}^{D}, \mathbf{C})$,
$A(t,x)$ is given as
\begin{eqnarray*}
A(t,x) = \frac{1}{(2\pi)^D} \int\! v_0(\hat{y}_1+\hat{x}_1) 
    \int\! e^{-i\hat{y}_1 \hat{\xi}_1} e^{Q(i\hat{\xi}_1, 0)t} d\hat{\xi}_1d\hat{y}_1.
\end{eqnarray*}
The rest of the proof is the same as those of Prop.3.1 and 3.2. \hfill $\blacksquare$

%%%%%%%%%%%%%%%%%%%%%%%%%%%%%%%%%%%%%%%%%%%%%%%%%%%%%%%%%%%%%%%%%%%%%%%%%%%%%%%%%%%%%%%%%%%%%%%%%%

\subsection{Higher dimensional case}

Suppose $u = (u_1, \cdots ,u_m) \in \mathbf{C}^m$ and $x = (x_1, \cdots ,x_d) \in \mathbf{R}^d$.
For fixed $1\leq D\leq d$, we use the same notation $x = (\hat{x}_1, \hat{x}_2)$ as in Sec.3.1.
Let $\{ P_{ij}(x) \}_{i,j=1}^m$ be the set of 
polynomials of $x$.
Define the matrix $P(x)$ by
\begin{equation}
P(x) = P(x_1 , \cdots ,x_d) = \left(
\begin{array}{@{\,}ccc@{\,}}
P_{11}(x) & \cdots & P_{1m}(x) \\
\vdots & \ddots & \vdots \\ 
P_{m1}(x) & \cdots & P_{mm}(x) \\
\end{array}
\right).
\end{equation}
The differential operator $\mathcal{P}$ is defined to be $\mathcal{P} = P(\partial _1, \cdots , \partial _d)$.
The algebraic equation
\begin{equation}
\det (\lambda - P(i\xi)) = \det \left(
\begin{array}{@{\,}ccc@{\,}}
\lambda - P_{11}(i\xi) & \cdots & -P_{1m}(i\xi) \\
\vdots & \ddots & \vdots \\ 
-P_{m1}(i\xi) & \cdots & \lambda -P_{mm}(i\xi) \\
\end{array}
\right) = 0
\end{equation}
is called the dispersion relation.
Let $\lambda _1(\xi) , \cdots ,\lambda _m(\xi)$ be roots of this equation.
Then, $\lambda _1(\mathbf{R})\cup \cdots  \cup \lambda _m(\mathbf{R})$ gives the spectrum of $\mathcal{P}$.
We suppose for simplicity that only $\lambda _1(\xi)$ contributes to the center subspace of $\mathcal{P}$
(see (C1) below).
Extending to more general situations is not difficult (see Remark 3.7 below).
\\[0.2cm]
\textbf{(C0)} The matrix $P(i\xi)$ is diagonalizable for any $\xi \in \mathbf{R}^d$.
\\
\textbf{(C1)} $\mathrm{Re}[\lambda _1(\xi)] \leq 0$ and $\mathrm{Re}[\lambda _j(\xi)] < 0$
for any $\xi \in \mathbf{R}^d$ and $j=2,\cdots , m$.
\\
\textbf{(C2)}\, There exist $\omega \in \mathbf{R},\, k\in \mathbf{R}^d$ and an integer $M$ such that
\begin{eqnarray*}
& & \lambda _1(k) = i\omega, \\
& & \frac{\partial ^\beta \lambda_1 }{\partial \hat{x}_1^\beta }(k) = 0, \,\,\text{for any $\beta $ such that}
\,\, |\beta | = 1, \cdots  ,M-1, \\
& & \frac{\partial ^\beta \lambda _1}{\partial \hat{x}_1^\beta }(k) \neq 0,\,\, \text{for some $\beta $ such that}
\,\, |\beta | = M.
\end{eqnarray*}
\textbf{(C3)} Define $Q(x)$ and $\mathcal{Q}$ by
\begin{equation}
Q(x) = Q(\hat{x}_1, 0) = \sum_{|\beta | = M} \frac{1}{(\beta _1!)\cdots (\beta _D !)} 
\frac{\partial ^\beta \lambda _1}{\partial \hat{x}_1^\beta }(k) (\hat{x}_1 / i)^\beta , \quad 
\mathcal{Q} = Q(\partial _1 , \cdots ,\partial _D, 0 ,\cdots ,0).
\end{equation}
Then, both of $\mathcal{P}$ and $\mathcal{Q}$ are elliptic in the sense that
there exist $c_1,c_2 > 0$ such that $\mathrm{Re}[\lambda _j(\xi )] < -c_2|\xi|^2 \,( j=1,\cdots ,m)$  and 
$\mathrm{Re}[Q(i\hat{\xi}_1)] < -c_2|\hat{\xi}_1|^2$ hold for $|\xi|, |\hat{\xi}_1| \geq c_1$,
where $\hat{\xi}_1 = (\xi_1, \cdots ,\xi_D)$.
\\

Put $B^r = BC^r(\mathbf{R}^r; \mathbf{C})$ and let $(B^r)^m$ be the product space.
The norm on $(B^r)^m$ is defined by $|| u || = \max_{1\leq j\leq m}|| u_j ||$
for $u=(u_1, \cdots ,u_m)$.
Note that $\mathcal{P}$ is an operator densely defined on $(B^r)^m$ while 
$\mathcal{Q}$ is an operator densely defined on $B^r$.
When $m=1$, $\lambda _1(\xi) = P(i\xi)$, so that the above assumptions and $Q$ are reduced to those in Sec.3.1.
\\[0.2cm]
\textbf{Example \thedef.}
Suppose $m=d=2$ and consider the operator $L$ defined by (\ref{1-15}) with the 
condition $0<d<k<1,\, (k+d)^2 = 4d$.
The dispersion relation is given by (\ref{2-8}), whose roots are denoted as $\lambda _2(\xi)< \lambda _1(\xi)$.
It is easy to verify that $\lambda _1(\xi) = 0$ if and only if $\xi = (\xi_1, \xi_2)$ satisfies 
$\xi_1^2 + \xi_2^2 = c^2 := (k-d)/2d$.
Thus there are infinitely many points $\xi$ satisfying $\lambda _1(\xi) = 0$.
We choose $(\xi_1, \xi_2) = (0,c)$.
Then, we can show that
\begin{eqnarray*}
\lambda _1(0,c) = \frac{\partial \lambda _1}{\partial \xi_1}(0,c)
= \frac{\partial^2 \lambda _1}{\partial \xi_1^2}(0,c) = \frac{\partial^3 \lambda _1}{\partial \xi_1^3}(0,c)=0,
\quad  \frac{\partial^4 \lambda _1}{\partial \xi_1^4}(0,c) \neq 0,
\end{eqnarray*}
while
\begin{eqnarray*}
\frac{\partial \lambda _1}{\partial \xi_2}(0,c) = 0, \quad 
\frac{\partial^2 \lambda _1}{\partial \xi_2^2}(0,c) \neq 0.
\end{eqnarray*}
Hence, (C2) is satisfied with $k=(0,c),\, \omega =0,\, D=1$ (i.e. $\hat{x}_1=x$ and $\hat{x}_2=y$), and $M=4$.
In this case, 
\begin{equation}
\mathcal{Q} = \frac{1}{4!} \frac{\partial^4 \lambda _1}{\partial \xi_1^4}(0,c) \frac{\partial ^4}{\partial x^4}
 = - \frac{2d^2}{(k+d)(1-d)}\frac{\partial ^4}{\partial x^4},
\label{3-17}
\end{equation}
see Eq.(\ref{1-17}).

Let $\bm{w} = (w_1, \cdots ,w_m)$ be an eigenvector of $P(ik)$ associated with $\lambda _1(k) = i\omega $.
Note that $e^{ikx}\bm{w}$ is an eigenfunction of $\mathcal{P}$ included in the center subspace.
Consider two systems of PDEs:
\begin{subnumcases}
{}
\displaystyle \frac{\partial u}{\partial t}  = \mathcal{P}u, \quad
\displaystyle u(0,x) = v_0(\hat{x}_1)e^{ikx}\bm{w}, \label{3-17a}\\[0.2cm]
\displaystyle \frac{\partial A}{\partial t} = \mathcal{Q}A. \quad
\displaystyle A(0,x) = v_0(\hat{x}_1). \label{3-17b}
\end{subnumcases}
Note that $v_0$ depends only on $\hat{x}_1 = (x_1, \cdots  ,x_D)$.
The former is a system of $m$-equations on $\mathbf{R}^d$, while the latter
is a single equation on $\mathbf{R}^D$.
We also consider the perturbative problem
\begin{subnumcases}
{}
\displaystyle \frac{\partial u}{\partial t}  = \mathcal{P}u , \quad
\displaystyle u(0,x) = v_0(\eta \hat{x}_1)e^{ikx}\bm{w}, \label{3-18a}\\[0.2cm]
\displaystyle \frac{\partial A}{\partial t} = \mathcal{Q}A, \quad
\displaystyle A(0,x) = v_0(\eta \hat{x}_1),\label{3-18b}
\end{subnumcases}
where $\eta = \varepsilon ^{1/M}$ and $\varepsilon >0$ is a small parameter.
Solutions of them satisfy the next proposition.
\\[0.2cm]
\textbf{Proposition \thedef.}
Suppose $r\geq 0$.
Under the assumptions (C0) to (C3), there exists a constant $C_1>0$ such that
\begin{equation}
|| e^{\mathcal{P}t}(e^{ikx} v_0\cdot \bm{w}) - e^{i\omega t + ikx}(e^{\mathcal{Q}t}v_0) \cdot \bm{w} ||
 \leq C_1 t^{-1/M} || v_0 ||
\label{3-20}
\end{equation}
holds for any $t>0$ and $v_0 \in BC^r(\mathbf{R}^r; \mathbf{C})$.
Next, suppose $r\geq 1$.
For any $\varepsilon >0$ and $t_0 > 0$, there exists a positive number $C_1 = C_1(t_0)$ such that the inequality
\begin{equation}
|| e^{\mathcal{P}t}(e^{ikx} \hat{v}_0\cdot \bm{w}) - e^{i\omega t + ikx}(e^{\mathcal{Q}t}\hat{v}_0)\cdot \bm{w}|| 
\leq \eta C_1 || v_0 ||
\label{3-21}
\end{equation}
holds for $t \geq t_0$ and $v_0 \in BC^r(\mathbf{R}^r; \mathbf{C})$, where $\hat{v}_0(\cdot) := v_0(\eta \,\cdot)$.
\\[0.2cm]
\textbf{Proof.}
We suppose $D = d$ for simplicity; that is, $\hat{x}_1 = x$ and $\beta = \alpha $.
The case $D < d$ is easily reduced to the case $D = d$ as in the proof of Prop.3.3.
We also suppose $\omega  = 0$ without loss of generality.

Like as the proof of Prop.3.1, a solution of (\ref{3-17a}) is written as 
\begin{eqnarray*}
u(t,x)= \frac{e^{ikx}}{(2\pi)^d} \int\! v_0(y+x) \int\! e^{-iy\xi} e^{P(i\xi + ik)t}\bm{w} d\xi dy.
\end{eqnarray*}
Note that $e^{P(i\xi + ik)t}$ is an exponential of a matrix.
Let $S(\xi)$ be a matrix such that
\begin{equation}
S(\xi)^{-1} P(i\xi) S(\xi) := \Lambda (\xi) = \left(
\begin{array}{@{\,}ccc@{\,}}
\lambda _1(\xi) & &  \\
& \ddots & \\
& & \lambda _m(\xi) 
\end{array}
\right) .
\label{3-22}
\end{equation}
Because of the assumption (C0), we can assume that $S(\xi), S(\xi)^{-1}$ and $\lambda _j(\xi)$'s
are smooth in $\xi$.
Then, 
\begin{eqnarray*}
u(t,x)= \frac{e^{ikx}}{(2\pi)^d} \int\! v_0(y+x) \int\! e^{-iy\xi} 
S(k+\xi) e^{\Lambda (k + \xi)t}S(k+\xi)^{-1} \bm{w} d\xi dy.
\end{eqnarray*}
Put $\tau = t^{-1/M}$.
Changing variables $\xi \mapsto \tau \xi,\, y\mapsto y/\tau$ yields
\begin{eqnarray*}
u(t,x)= \frac{e^{ikx}}{(2\pi)^d} \int\! v_0(y/\tau +x) \int\! e^{-iy\xi} 
S(k + \tau \xi) e^{\Lambda (k + \tau \xi )/\tau^M}S(k+\tau \xi )^{-1} \bm{w} d\xi dy.
\end{eqnarray*}
Expanding $S(k + \tau \xi)$ and $S(k + \tau \xi)^{-1}$, it turns out that there is a function $G_1$
such that
\begin{eqnarray*}
u(t,x)= \! \frac{e^{ikx}}{(2\pi)^d} \!\!\int\!\! v_0(y/\tau +x) \!\!\int\!\! e^{-iy\xi} 
S(k) e^{\Lambda (k + \tau \xi )/\tau^M}S(k)^{-1} \bm{w} d\xi dy + 
\tau \!\int\!\! v_0(y/\tau +x) G_1(y, \tau) dy.
\end{eqnarray*}
By a similar estimate used in the proof of Prop.3.1,
we can show that there exists $D_1>0$ such that the norm of the second term above 
has an upper bound $ \tau D_1 || v_0 ||$ for any $\tau > 0$.
Since $\bm{w}$ is an eigenvector associated with $\lambda _1(\xi)$, we obtain
\begin{eqnarray*}
u(t,x)= \frac{e^{ikx}}{(2\pi)^d} \int\! v_0(y/\tau +x) \int\! e^{-iy\xi} 
 e^{\lambda_1 (k + \tau \xi )/\tau^M} \bm{w} d\xi dy + \tau \!\int\!\! v_0(y/\tau +x) G_1(y, \tau) dy.
\end{eqnarray*}
Therefore, we have
\begin{eqnarray*}
& & u(t,x) - e^{ikx}A(t,x)\bm{w}  \\
&=& 
\frac{e^{ikx}}{(2\pi)^d} \int\! v_0(y/\tau +x) \int\! e^{-iy\xi} e^{Q(i\xi)} 
\left( e^{g(\xi, \tau)}  -1   \right) \bm{w} d\xi dy 
+ \tau \!\int\!\! v_0(y/\tau +x) G_1(y, \tau) dy,
\end{eqnarray*}
where $g(\xi, \tau)  = \lambda_1 (k + \tau \xi )/\tau^M - Q(i\xi)$.
The rest of the proof is the same as those of Prop.3.1 and 3.2. \hfill $\blacksquare$
\\[0.2cm]
\textbf{Remark \thedef.} Let $\lambda _1 (\xi), \cdots ,\lambda _m(\xi)$ be eigenvalues of $P(i\xi)$ as before.
Even if several eigenvalues lie on the imaginary axis and (C1) is violated,
to modify Prop.3.6 is very easy;
since equations are linear, the superposition principle is applicable.
A typical problem is that $P(i\xi)$ is real-valued and eigenvalues occur in complex conjugate pairs.
For example, suppose (C0), (C2), (C3) and the following (C1)' instead of (C1):
\\[0.2cm]
\textbf{(C1)'} $\lambda _2(\xi) = \overline{\lambda _1(\xi)}$. 
$\mathrm{Re}[\lambda _{1,2}(\xi)] \leq 0$ and $\mathrm{Re}[\lambda _j(\xi)] < 0$
for any $\xi \in \mathbf{R}^d$ and $j=3,\cdots , m$.
\\[0.2cm]
Put $\overline{\mathcal{Q}} : = \overline{Q(\partial _1, \cdots ,\partial _D, 0, \cdots ,0)}$.
Let $\bm{w}_1$ and $\bm{w}_2$ be eigenvectors of $P(ik)$ associated with $\lambda _1(k) = i\omega $ and $\lambda _2(k) = -i\omega $,
respectively.
In this case, instead of Eq.(\ref{3-20}), the inequality
\begin{equation}
|| e^{\mathcal{P}t}(e^{ikx} v_1\cdot \bm{w}_1 + e^{ikx} v_2\cdot \bm{w}_2)
  - e^{i\omega t + ikx}(e^{\mathcal{Q}t}v_1) \cdot \bm{w}_1  - e^{-i\omega t + ikx}(e^{\overline{\mathcal{Q}}t}v_2) \cdot \bm{w}_2||
 \leq C_1 t^{-1/M} (|| v_1 || + || v_2 ||)
\label{3-23}
\end{equation}
holds for any $t>0$ and $v_1, v_2 \in BC^r(\mathbf{R}^r; \mathbf{C})$, and similarly for Eq.(\ref{3-21}).
\\[0.2cm]
\textbf{Example \thedef.}
Suppose $m=2$ and $d=1$.
Define a linear operator
\begin{equation}
\mathcal{P}\left(
\begin{array}{@{\,}c@{\,}}
u \\
v
\end{array}
\right) = \left(
\begin{array}{@{\,}cc@{\,}}
D \partial ^2 & 1 \\
-1 & D \partial ^2
\end{array}
\right) \left(
\begin{array}{@{\,}c@{\,}}
u \\
v
\end{array}
\right),
\end{equation}
where $D > 0$ is a diffusion constant.
This operator arises from Eq.(\ref{1-7b}).
Eigenvalues of $P(i\xi)$ are $\lambda _1 (\xi) = -D \xi^2 + i$ and $\lambda _2 (\xi) = -D \xi^2 - i$.
Hence, the assumptions (C0), (C1)', (C2) and (C3) are satisfied with
\begin{eqnarray*}
& & \lambda _1(0) = i, \,\, \lambda _2(0) = -i,\,\, M = 2,\,\, \mathcal{Q} = \overline{\mathcal{Q}} = D \frac{\partial ^2}{\partial x^2}, \\
& & \bm{w}_1 = \left(
\begin{array}{@{\,}c@{\,}}
1 \\
i
\end{array}
\right), \quad \bm{w}_2 = \left(
\begin{array}{@{\,}c@{\,}}
1 \\
-i
\end{array}
\right) = \overline{\bm{w}}_1.
\end{eqnarray*}
Eq.(\ref{3-23}) is given as
\begin{equation}
|| e^{\mathcal{P}t}(v_1\cdot \bm{w}_1 + v_2\cdot \bm{w}_2)
  - e^{i\omega t}(e^{\mathcal{Q}t}v_1) \cdot \bm{w}_1  - e^{-i\omega t}(e^{\mathcal{Q}t}v_2) \cdot \bm{w}_2||
 \leq C_1 t^{-1/M} (|| v_1 || + || v_2 ||).
\end{equation}
In most applications, we take $v_2(x) = \overline{v_1(x)}$ to obtain a real-valued solution of 
$\dot{\bm{u}} = \mathcal{P}\bm{u}$.
The above inequality implies that an approximate solution of $\dot{\bm{u}} = \mathcal{P}\bm{u}$
is constructed through the complex heat equation $\dot{A} = \mathcal{Q}A$.

%%%%%%%%%%%%%%%%%%%%%%%%%%%%%%%%%%%%%%%%%%%%%%%%%%%%%%%%%%%%%%%%%%%%%%%%%%%%%%%%%%%%%%%%%%%%%%%%%%
%%%%%%%%%%%%%%%%%%%%%%%%%%%%%%%%%%%%%%%%%%%%%%%%%%%%%%%%%%%%%%%%%%%%%%%%%%%%%%%%%%%%%%%%%%%%%%%%%%

\section{Main theorems}

In this section, a reduction of a perturbation term $f(u)$ is given.
Combined with the reduction of linear semigroups, a reduction of Eq.(\ref{1-1}) is performed.

%%%%%%%%%%%%%%%%%%%%%%%%%%%%%%%%%%%%%%%%%%%%%%%%%%%%%%%%%%%%%%%%%%%%%%%%%%%%%%%%%%%%%%%%%%%%%%%%%%

\subsection{One dimensional case}

We start with the case $u\in \mathbf{C}$ and $x\in \mathbf{R}^d$.
Put $x = (x_1, \cdots ,x_d)$ and $\alpha = (\alpha _1, \cdots ,\alpha _d)$, 
where $\alpha$ denotes a multi-index.
For a fixed integer $1\leq D\leq d$, we denote $x \in \mathbf{R}^d$ as 
$x = (\hat{x}_1, \hat{x}_2)$ with $\hat{x}_1 = (x_1, \cdots , x_D)$
and $\hat{x}_2 = (x_{D+1} , \cdots ,x_d)$.
Accordingly, a multi-index $\alpha $ is also denoted as $\alpha  = (\beta, \gamma )$.
Let $P(x) = \sum^q_{|\alpha |=0} a_\alpha x^\alpha $ be a polynomial and 
$\mathcal{P}:= P(\partial _1 , \cdots ,\partial _d)$ a differential operator on $\mathbf{R}^d$,
where $\partial _j$ denotes the derivative with respect to $x_j$.
For the main theorems, we make the following assumptions.
\\[0.2cm]
\textbf{(D1)} $\mathrm{Re}[P(i\xi)] \leq 0$ for any $\xi \in \mathbf{R}^d$.
\\
\textbf{(D2)}\, There exist $\omega \in \mathbf{R},\, k\in \mathbf{R}^d \,\, ( (\omega , k) \neq (0,0))$, 
a finite set of integers $J = \{ j_1, \cdots , j_N \}$ and $\{M_1, \cdots ,M_N\}$ such that
\begin{eqnarray*}
& & P(ij_nk) = ij_n\omega,\,\, (n=1,\cdots ,N), \\
& & \frac{\partial ^\beta P}{\partial \hat{x}_1^\beta }(ij_nk) = 0, \,\,\text{for any $\beta $ such that}
\,\, |\beta | = 1, \cdots  ,M_n-1, \,\, (n=1,\cdots ,N),\\
& & \frac{\partial ^\beta P}{\partial \hat{x}_1^\beta }(ij_nk) \neq 0,\,\, \text{for some $\beta _n$ such that}
\,\, |\beta _n| = M_n, \,\, (n=1,\cdots ,N).
\end{eqnarray*}
The set $J$ consists of \textit{all} integers satisfying $P(ijk) = ij\omega$.
\\
\textbf{(D3)} For $n=1,\cdots ,N$, define $Q_n(x)$ and $\mathcal{Q}_n$ by
\begin{equation}
Q_n(x) = Q_n(\hat{x}_1, 0) = \sum_{|\beta | = M_n} \frac{1}{(\beta _1!)\cdots (\beta _D !)} 
\frac{\partial ^\beta P}{\partial \hat{x}_1^\beta }(ij_nk) \hat{x}_1^\beta , \quad 
\mathcal{Q}_n = Q_n(\partial _1 , \cdots ,\partial _D, 0 ,\cdots ,0).
\end{equation}
Then, both of $\mathcal{P}$ and $\mathcal{Q}_n$ are elliptic in the sense that there exist $c_1,c_2>0$ such that
$\mathrm{Re}[P(i\xi )] < -c_2|\xi|^2$ and $\mathrm{Re}[Q_n(i\hat{\xi}_1,0)] < -c_2|\hat{\xi}_1|^2$ hold for 
$|\xi|, |\hat{\xi}_1| \geq c_1$ and $n=1,\cdots ,N$, where $\hat{\xi}_1 = (\xi_1, \cdots ,\xi_D)$.
\\

In addition to (B2) before, a new assumption $(\omega ,k) \neq (0,0)$ and the set $J$ are introduced.
In most examples, $J$ consists of $J = \{ +1, -1\}$ as Sec.1, see also an example below.

Put $B^r = BC^r(\mathbf{R}^d; \mathbf{C})$.
For a given function $f : B^r\to B^r$, let us consider the Fourier series of the quantity
$f(\sum^N_{n=1} A_n e^{ij_n \omega t + ij_n kx})$, where $A = (A_1, \cdots  ,A_N ) \in \mathbf{C}^N$.
Since $(\omega ,k) \neq (0,0)$, the Fourier series is well-defined and it is easy to verify that 
the series is of the form
\begin{equation}
f(\sum^N_{n=1} A_n e^{ij_n \omega t + ij_n kx})=\sum^\infty_{j=-\infty}C_j(A)e^{ij\omega t + ijkx}.
\label{4-3}
\end{equation}
For example, when $k_1 \neq 0$, $C_j(A)$ is given by
\begin{equation}
C_j(A) = C_j(A_1, \cdots ,A_N) := 
\frac{k_1}{2\pi} \int^{2\pi / k_1}_{0}\! f(\sum^N_{n=1} A_n e^{ij_n k_1x_1}) e^{-ijk_1x_1}dx_1.
\label{4-2}
\end{equation}
When $\omega \neq 0$, it is also written as
\begin{equation}
C_j(A) =
\frac{\omega }{2\pi} \int^{2\pi / \omega }_{0}\! f(\sum^N_{n=1} A_n e^{ij_n \omega t}) e^{-ij \omega t}dt.
\end{equation}
In particular, $C_{j_n} (A)$ is denoted by $R_n(A)$ if $j_n\in J$.
For any $\theta \in \mathbf{R}$ and $j\in \mathbf{Z}$, $C_j(A)$ satisfies the equality
\begin{equation}
C_j(e^{ij_1 \theta }A_1 ,\cdots ,e^{ij_N\theta }A_N) = e^{ij\theta }C_j(A_1 ,\cdots ,A_N).
\label{4-3b}
\end{equation}
Let $\varepsilon >0$ be a small parameter.
We will consider the two initial value problems:
\begin{subnumcases}
{}
\displaystyle \frac{\partial u}{\partial t}  = \mathcal{P}u + \varepsilon f(u) , \quad
\displaystyle u(0,x) = \sum^N_{n=1} e^{ij_n kx}v_n(\eta \hat{x}_1), \label{4-4a}\\[0.2cm]
\displaystyle \frac{\partial A_n}{\partial t} = \mathcal{Q}_nA_n + \varepsilon R_n(A), \quad
\displaystyle A_n(0,x) = v_n(\eta \hat{x}_1),\,\, (n= 1,\cdots ,N), \label{4-4b}
\end{subnumcases}
where $\eta = \varepsilon ^{1/M}$ and $M := \min \{ M_1, \cdots ,M_N\}$.
Note that the former is a single equation while the latter is a system of PDEs.
\\[0.2cm]
\textbf{Example \thedef.}
Let us consider the Swift-Hohenberg equation (\ref{1-2}).
For this equation, $D = d= 1$ and $P(x) = -(x^2 + k^2)^2$.
Since $P(i\xi) = 0$ if and only if $\xi = \pm k$, the set $J$ consists of $j_1 = +1,\, j_2 = -1$.
We have
\begin{eqnarray*}
\frac{\partial P}{\partial x}(\pm ik) = 0, \quad \frac{\partial^2 P}{\partial x^2}(\pm ik) = 8k^2.
\end{eqnarray*} 
Thus $M_1 = M_2 = 2$, and both of $\mathcal{Q}_1$ and $\mathcal{Q}_2$ are given by
\begin{eqnarray*}
\mathcal{Q}_{1,2} = \frac{1}{2}\frac{\partial^2 P}{\partial x^2}(\pm ik)\frac{\partial ^2}{\partial x^2}
 = 4k^2 \frac{\partial ^2}{\partial x^2}.
\end{eqnarray*}
Since $f(u) = u-u^3$, the expansion of $f(A_1 e^{ikx} + A_2e^{-ikx})$ is
\begin{eqnarray*}
f(A_1 e^{ikx} + A_2e^{-ikx})
 = A_1e^{ikx} + A_2e^{-ikx} - (A_1^3 e^{3ikx} + 3A_1^2A_2 e^{ikx} + 3A_1A_2^2 e^{-ikx}+A_2^3 e^{-3ikx} ).
\end{eqnarray*}
This provides
\begin{eqnarray*}
R_1(A) = C_1(A) = A_1 - 3A_1^2A_2, \quad R_2(A) = C_{-1}(A) = A_2 - 3A_1A_2^2.
\end{eqnarray*}
Therefore, the amplitude equation (\ref{4-4b}) is given by
\begin{equation}
\frac{\partial A_1}{\partial t} = 4k^2 \frac{\partial ^2A_1}{\partial x^2} + \varepsilon (A_1 - 3A_1^2A_2),
\quad 
\frac{\partial A_2}{\partial t} = 4k^2 \frac{\partial ^2A_2}{\partial x^2} + \varepsilon (A_2 - 3A_1A_2^2).
\label{4-5}
\end{equation}
Usually, we assume $A_1 = \overline{A_2}$, which gives the Ginzburg-Landau equation (\ref{1-3}).
\\

Put $\hat{v}_n(x) = v_n (\eta x)$.
The above equations are rewritten as integral equations of the form
\begin{subnumcases}
{}
\displaystyle u = e^{\mathcal{P}t} (\sum^N_{n=1}e^{ij_nkx}\hat{v}_n)
        + \varepsilon \int^{t}_{0}\! e^{\mathcal{P}(t-s)} f(u(s))ds , \label{4-6a}\\[0.2cm]
\displaystyle A_n = e^{\mathcal{Q}_nt}\hat{v}_n 
        + \varepsilon \int^{t}_{0}\! e^{\mathcal{Q}_n(t-s)}R_n(A(s))ds,\,\,  (n= 1,\cdots ,N), \label{4-6b}
\end{subnumcases}
whose solutions are called mild solutions.
When $f : B^r\to B^r$ is $C^1$, then $R_n : B^r\to B^r$ is also $C^1$, and
due to the standard existence theorem (see Pazy\cite{Paz}), there exists a positive number $T_0 > 0$
such that the above integral equations have mild solutions $u(t, \cdot),\, A_n(t, \cdot) \in B^r$
for $0\leq t \leq T_0/\varepsilon $.
In particular, when the initial condition $\{ v_n\}$ is included in the domain of $\mathcal{P}$
and $\mathcal{Q}_n$, then a mild solution is a classical solution which is differentiable in $t>0$.
In this paper, we only consider mild solutions.
The main theorems for a one-dimensional case are stated as follows:
\\[0.2cm]
\textbf{Theorem \thedef.} 
Suppose $f : BC^r(\mathbf{R}^d; \mathbf{C}) \to BC^r(\mathbf{R}^d; \mathbf{C})\,\, (r\geq 1)$ is $C^1$ and $\varepsilon >0$ is sufficiently small.
For any $\{ v_n\}^N_{n=1} \subset BC^r(\mathbf{R}^d; \mathbf{C})$,
there exist positive numbers $C,T_0$ and $t_0$ such that mild solutions of 
the two initial value problems (4.5) satisfy
\begin{equation}
|| u(t,x) - \sum^N_{n=1}A_n(t,x) e^{ij_n\omega t + ij_nkx} || 
\leq C\eta = C\varepsilon ^{1/M},
\end{equation}
for $t_0 \leq t \leq T_0/\varepsilon $.
\\

Next, let us show that the error estimate above holds for any $t>t_0$ under a suitable condition.
For ordinary differential equations, it is proved in \cite{Chi} that if the amplitude equation has a stable
hyperbolic invariant manifold, then a given equation has a stable invariant manifold of the same type and 
approximate solutions are valid for any $t>0$ near the manifold.
For our situation, suppose that there is a constant vector $\phi \in \mathbf{R}^N$ such that $R(\phi) = 0
,\, R = (R_1, \cdots ,R_N)$.
Then, $\phi$ is a steady state of the amplitude equation.
Unfortunately, $\phi$ is not hyperbolic because of the symmetry (\ref{4-3b});
the Jacobi matrix of $R$ at $\phi$ has a zero-eigenvalue in general.
For example, although the amplitude equation (\ref{4-5}) for the Swift-Hohenberg equation
has a steady state $(A_1, A_2) = (1/\sqrt{3}, 1/\sqrt{3})$, 
the Jacobi matrix of $R$ at $(1/\sqrt{3}, 1/\sqrt{3})$ has a zero-eigenvalue.
However, if we restrict solutions to the invariant set $\{ A_1 = A_2\}$, (\ref{4-5}) is reduced to
\begin{equation}
\frac{\partial A}{\partial t} = 4k^2 \frac{\partial ^2A}{\partial x^2} + \varepsilon (A-3A^3),
\label{4-37}
\end{equation}
and the derivative of the function $A-3A^3$ at $A = 1/\sqrt{3}$ is negative.
This implies that $A = 1/\sqrt{3}$ is a hyperbolically stable steady state of (\ref{4-37}),
and we expect that the  Swift-Hohenberg equation also has a corresponding stable solution.
For more general situations, we make the following assumption.
\\[0.2cm]
\textbf{(D4)} For $|\beta| = M_n$ and $ n=1,\cdots ,N$,
\[ P(i\xi) = \overline{P(-i\xi)},\,\, f(\overline{u}) = \overline{f(u)}\quad 
\text{and}  \quad \displaystyle \frac{\partial^\beta  P}{\partial \hat{x}_1^\beta}(ij_nk) = 
\frac{\partial^\beta  P}{\partial \hat{x}_1^\beta}(-ij_nk).\]

The first two equalities imply that $P(x)$ and $f(u)$ are real-valued when $x, u\in \mathbf{R}$.
Due to this assumption, $P(-ij_n k) = -ij_n\omega $ when $P(ij_n k) = ij_n\omega $.
Hence, the set $J$ of integers satisfying $P(ijk) = ij\omega $ is given by
$J = \{ j_1, \cdots ,j_N\} \cup \{ -j_1, \cdots ,-j_N\}$.
We denote $-j_n$ by $j_{-n}$.
Then, $M_n$ and $\mathcal{Q}_n$ are defined for $n = \pm 1,\cdots ,\pm N$ as in (D2), (D3).
It follows from (D4) that $M_n = M_{-n}$ and $\mathcal{Q}_n = \mathcal{Q}_{-n}$.
For many examples, $J$ consists of $J=\{ \pm 1\}$ and this assumption is satisfied.
In the present notation, the function $C_j(A)$ defined by (\ref{4-2}) is given by
\begin{eqnarray}
\!\!\!\!\!C_j(A_1, \cdots ,A_N, A_{-1}, \cdots ,A_{-N}) \!
= \frac{k_1}{2\pi} \!\!\int^{2\pi / k_1}_{0}\!\!\! f(\sum^N_{n=1} A_n e^{ij_n k_1x_1} 
\!+\! \sum^N_{n=1}A_{-n}e^{-ij_n k_1x_1})  e^{-ijk_1x_1} dx_1.
\label{4-35b}
\end{eqnarray}
In particular, $C_{j_n}(A)$ is denoted by $R_n(A)$ for $n = \pm 1,\cdots ,\pm N$.
Hence, the amplitude equation is given as a system of $2N$-equations of the form
\begin{equation}
\frac{\partial A_n}{\partial t} = \mathcal{Q}_n A_n + \varepsilon R_n(A), \quad (n = \pm 1,\cdots ,\pm N ).
\label{4-35}
\end{equation}
This system can be reduced as follows; It is easy to verify that $C_j$ satisfies
\begin{equation}
C_j(A_1, \cdots ,A_N, A_{-1}, \cdots ,A_{-N}) = C_{-j}(A_{-1}, \cdots ,A_{-N}, A_1, \cdots ,A_N).
\label{4-35c}
\end{equation}
Thus putting $A_n = A_{-n}$ yields $R_n(A) = R_{-n}(A)$.
Since $\mathcal{Q}_n = \mathcal{Q}_{-n}$, putting $A_n = A_{-n}$ shows that Eq.(\ref{4-35})
is reduced to the system of $N$-equations
\begin{equation}
\frac{\partial A_n}{\partial t} = \mathcal{Q}_n A_n + 
       \varepsilon R_n(A_1, \cdots ,A_N, A_1, \cdots ,A_N), \quad (n = 1,\cdots ,N).
\label{4-36}
\end{equation}

Define the function $S_n$ to be
\begin{equation}
S_n(A) = S_n(A_1, \cdots ,A_N) = R_n(A_1, \cdots ,A_N, A_1, \cdots ,A_N).
\label{4-13b}
\end{equation}
We consider the two initial value problems:
\begin{subnumcases}
{}
\displaystyle \frac{\partial u}{\partial t}  = \mathcal{P}u + \varepsilon f(u) , \quad
\displaystyle u(0,x) = \sum^N_{n=1} \left( e^{ij_n kx}
                     +  e^{-ij_n kx}\right) v_n(\eta \hat{x}_1), \label{4-39a}\\[0.2cm]
\displaystyle \frac{\partial A_n}{\partial t} = \mathcal{Q}_nA_n + \varepsilon S_n(A), \quad
\displaystyle A_n(0,x) = v_n(\eta \hat{x}_1),\,\, (n= 1,\cdots ,N), \label{4-39b}
\end{subnumcases}
where $\eta = \varepsilon ^{1/M}$ and $M := \min \{ M_1, \cdots ,M_N\}$.
Due to Thm.4.2, solutions of them satisfy
\begin{eqnarray*}
|| u(t,x) - \sum^N_{n=1}A_n(t,x)(e^{ij_n\omega t+ij_nkx} + e^{-ij_n\omega t-ij_nkx}) || \leq C\eta,
\end{eqnarray*}
for $t_0 \leq t\leq T_0/\varepsilon $.
Further, we can show the next theorem, in which $B^r = BC^r(\mathbf{R}^d;\mathbf{R})$ denotes the set of 
real-valued functions.
\\[0.2cm]
\textbf{Theorem \thedef.}  
Suppose (D1) to (D4) and 
$f :BC^r(\mathbf{R}^d;\mathbf{R}) \to BC^r(\mathbf{R}^d;\mathbf{R})\,\, (r\geq 1)$ is $C^2$ such that the second derivatives are locally Lipschitz continuous.
Suppose that there exists a constant vector $\phi = (\phi_1, \cdots ,\phi_N) \in \mathbf{R}^N$ such that
\\
\textbf{(i)} $S_n(\phi) = 0$ for $n=1, \cdots ,N$,
\\
\textbf{(ii)} the Jacobi matrix of $(S_1,\cdots ,S_N)$ at $\phi$ is diagonalizable and 
all eigenvalues of the matrix have negative real parts.
\\
If $\varepsilon >0$ is sufficiently small, Eq.(\ref{4-39a}) has a solution of the form
\begin{equation}
u_p(t,x) = \sum^N_{n=1} \Bigl( \phi_n + \eta \psi_n (t, x, \eta) \Bigr) \cdot 
 (e^{ij_n\omega t + ij_nkx} + e^{-ij_n\omega t - ij_nkx}).
\label{4-16c}
\end{equation}
The functions $\psi_n $ and $u_p$ are bounded as $\eta \to 0$ and satisfy
\begin{eqnarray*}
\left\{ \begin{array}{l}
\text{$2\pi / \omega $-periodic in $t$ (when $\omega \neq 0$)}, \\
\text{constant in $t$ (when $\omega = 0$)},  \\
\end{array} \right.
\quad
\left\{ \begin{array}{l}
\text{$2\pi / k_j $-periodic in $x_j$ (when $k_j \neq 0$)}, \\
\text{constant in $x_j$ (when $k_j = 0$)},  \\
\end{array} \right.
\end{eqnarray*}
for $j=1,\cdots ,d$.
This $u_p$ is stable in the following sense:
For any $n=1,\cdots ,N$, there is a neighborhood $U_n\subset BC^r(\mathbf{R}^d;\mathbf{R})$ of $\phi_n$ in 
$BC^r(\mathbf{R}^d;\mathbf{R})$ such that if $v_n\in U_n$,
then a mild solution $u$ of the initial value problem (\ref{4-39a}) satisfies
$|| u(t, \cdot) - u_p(t, \cdot) || \to 0$ as $t\to \infty$.
\\

The above conditions (i),(ii) show that $A(t,x) \equiv \phi$ is an asymptotically 
stable steady state of Eq.(\ref{4-39b}). 
Thus this theorem implies that a stable steady state of Eq.(\ref{4-39b}) induces a periodic solution of
Eq.(\ref{4-39a}).
Due to the symmetry (\ref{4-3b}), Eq.(\ref{4-35}) has a steady state $A_n = e^{ij_n \theta } \phi_n\,\,
(n=\pm 1,\cdots ,\pm N)$ for any $\theta \in [0, 2\pi)$.
Accordingly, we can prove that Eq.(\ref{4-39a}) has a stable periodic solution 
\begin{eqnarray*}
u_p(t,x) = \sum^N_{n=1} \Bigl( \phi_n + \eta \psi_n (t, x, \eta) \Bigr) \cdot 
 (e^{ij_n\omega t + ij_nkx + i j_n \theta } + e^{-ij_n\omega t - ij_nkx - ij_n \theta }),
\end{eqnarray*}
for any $\theta \in [0, 2\pi)$, the proof of which is reduced to (\ref{4-16c}) by the translation of $t$ or $x$.
\\[0.2cm]
\textbf{Proof of Thm.4.2.}
We prove the theorem for the case $D=d$ (thus $\hat{x}_1=x$ and $\alpha = \beta$) for simplicity of notation.
The general case $D < d$ can be proved in the same way.
A proof is divided into four steps.

\textit{Step 1. notation.}
It is convenient to introduce some notation:
We define a new coordinate $(T, X)$ by
\begin{eqnarray*}
& & x = X/\eta,\,\,  t = T/\varepsilon ,\,\, \hat{u}(T,X) = u(t,x),\,\,\hat{A}_n(T,X) = A_n(t,x), \\
& & \hat{\mathcal{P}} = \frac{1}{\varepsilon }P(\eta \partial _X),\,\,
\hat{\mathcal{Q}}_n = \frac{1}{\varepsilon }Q_n(\eta \partial _X).
\end{eqnarray*}
Then, Eqs.(\ref{4-4a}) and (\ref{4-4b}) are rewritten as
\begin{subnumcases}
{}
\displaystyle \frac{\partial \hat{u}}{\partial T}  = \hat{\mathcal{P}}\hat{u} +f(\hat{u}) , \quad
\displaystyle \hat{u}(0,X) = \sum^N_{n=1} e^{ij_n kX/\eta }v_n(X), \label{4-9a}\\[0.2cm]
\displaystyle \frac{\partial \hat{A}_n}{\partial T} = \hat{\mathcal{Q}}_n\hat{A}_n + R_n(\hat{A}), \quad
\displaystyle \hat{A}_n(0,X) = v_n(X),\,\, (n= 1,\cdots ,N). \label{4-9b}
\end{subnumcases}
Integrating them yields
\begin{subnumcases}
{}
\displaystyle \hat{u} = e^{\hat{\mathcal{P}}T} (\sum^N_{n=1}e^{ij_nkX/\eta }v_n)
        + \int^{T}_{0}\! e^{\hat{\mathcal{P}}(T-s)} f(\hat{u}(s))ds , \label{4-10a}\\[0.2cm]
\displaystyle \hat{A}_n = e^{\hat{\mathcal{Q}}_nT}v_n 
        + \int^{T}_{0}\! e^{\hat{\mathcal{Q}}_n(T-s)}R_n(\hat{A}(s))ds,\,\,  (n= 1,\cdots ,N), \label{4-10b}
\end{subnumcases}
which have mild solutions in $B^r$ for $0\leq T \leq T_0$.
For the function space $B^r$ written in the $X$-variable,
we introduce the norm $|| \cdot ||_\eta$ by
\begin{equation}
|| \varphi ||_\eta = \max_{0\leq |\alpha | \leq r} \sup_{X\in \mathbf{R}^d}
\{ \,\,\eta^{|\alpha |}\, |\partial ^\alpha \varphi (X)| \,\,\}.
\label{4-11}
\end{equation}
If we put $\hat{\varphi }(X) := \varphi (X/\eta ) = \varphi (x)$ for a given $\varphi (x)$, it is easy to see that 
$|| \varphi  || = || \hat{\varphi } ||_\eta$, where $|| \varphi  ||$ represents the standard norm on $B^r$.
In the present notation, Prop.3.2 is restated as follows:
there exist $t_0, C_1 > 0$ such that the inequality
\begin{equation}
|| e^{\hat{\mathcal{P}}T} (e^{ij_n kX/\eta} v_n) - e^{ij_n \omega T/\varepsilon + ij_n kX/\eta}
e^{\hat{\mathcal{Q}}_nT} v_n ||_\eta \leq \eta C_1 || v_n||_\eta
\label{4-12}
\end{equation}
holds for $T \geq \varepsilon t_0$ and for each $n= 1,\cdots ,N$.
Let us estimate $\hat{u} - \sum^N_{n=1} \hat{A}_n e^{ij_n\omega t+ ij_nkx}$
by using the norm $|| \cdot ||_\eta$.

\textit{Step 2. Gronwall inequality.}
It follows from Eq.(4.18) that
\begin{eqnarray*}
& & \hat{u} - \sum^N_{n=1} \hat{A}_n e^{ij_n\omega T/\varepsilon + ij_nkX/\eta } 
= \sum^N_{n=1}e^{\hat{\mathcal{P}}T} (e^{ij_nkX/\eta }v_n) 
         -\sum^N_{n=1}e^{ij_n\omega T/\varepsilon + ij_nkX/\eta }e^{\hat{\mathcal{Q}}_nT}v_n  \\
& &  \quad\quad        + \int^{T}_{0}\! e^{\hat{\mathcal{P}}(T-s)} f(\hat{u}(s))ds
       - \sum^N_{n=1} \int^{T}_{0}\!e^{ij_n\omega T/\varepsilon + ij_nkX/\eta } e^{\hat{\mathcal{Q}}_n(T-s)}R_n(\hat{A}(s))ds \\
&=& F(T) + \int^{T}_{0}\! e^{\hat{\mathcal{P}}(T-s)} 
   \Bigl(f(\hat{u}(s)) - f(\textstyle \sum^N_{n=1} \hat{A}(s)e^{ij_n\omega s/\varepsilon  + ij_nkX/\eta })\Bigr) ds,
\end{eqnarray*}
where 
\begin{eqnarray*}
F(T) &=& \sum^N_{n=1}e^{\hat{\mathcal{P}}T} (e^{ij_nkX/\eta }v_n) 
         -\sum^N_{n=1}e^{ij_n\omega T/\varepsilon + ij_nkX/\eta }e^{\hat{\mathcal{Q}}_nT}v_n \\
&+& \!\!\!\!\! \int^{T}_{0}\!\! e^{\hat{\mathcal{P}}(T-s)}
f(\textstyle \sum^N_{n=1} \hat{A}(s)e^{ij_n\omega s/\varepsilon  + ij_nkX/\eta }) ds \displaystyle
-\sum^N_{n=1} \int^{T}_{0}\!\!e^{ij_n\omega T/\varepsilon + ij_nkX/\eta } e^{\hat{\mathcal{Q}}_n(T-s)}R_n(\hat{A}(s))ds.
\end{eqnarray*}
Because of the existence theorem of mild solutions, there exists a positive constant $D_1$, 
which is independent of $\varepsilon $, such that 
\begin{eqnarray*}
|| \hat{u}(T) ||_\eta \leq D_1, 
\quad || \sum^N_{n=1} \hat{A}_n(T) e^{ij_n\omega T/\varepsilon + ij_nkX/\eta } ||_\eta \leq D_1,
\quad || e^{\hat{\mathcal{P}}T} ||_\eta \leq D_1
\end{eqnarray*}
hold for $0\leq T\leq T_0$.
Let $L > 0$ be a Lipschitz constant of $f$ in the ball $\{ \varphi \in B^r \, | \, || \varphi  ||_\eta \leq D_1\}$.
Then, we obtain
\begin{eqnarray*}
|| \hat{u} - \sum^N_{n=1} \hat{A}_n e^{ij_n\omega T/\varepsilon + ij_nkX/\eta }  ||_\eta \leq 
|| F(T) ||_\eta + \int^{T}_{0}\! D_1L  \,
|| \hat{u}(s) - \textstyle \sum^N_{n=1} \hat{A}_n(s)e^{ij_n\omega s/\varepsilon  + ij_nkX/\eta } ||_\eta ds.
\end{eqnarray*}
Gronwall inequality yields
\begin{eqnarray*}
|| \hat{u} - \sum^N_{n=1} \hat{A}_n e^{ij_n\omega T/\varepsilon + ij_nkX/\eta }  ||_\eta \leq 
|| F(T) ||_\eta + D_1L \int^{T}_{0}\! e^{D_1L(T-s)} || F(s)||_\eta ds.
\end{eqnarray*}
To estimate $F(T)$, we rewrite it with the aid of Eq.(\ref{4-3}) and $C_{j_n} = R_n$ as
\begin{eqnarray}
F(T) &=& \sum^N_{n=1}e^{\hat{\mathcal{P}}T} (e^{ij_nkX/\eta }v_n) 
         -\sum^N_{n=1}e^{ij_n\omega T/\varepsilon + ij_nkX/\eta }e^{\hat{\mathcal{Q}}_nT}v_n \nonumber \\
&+& \sum^N_{n=1} \int^{T}_{0}\! e^{ij_n \omega s/\varepsilon } 
    \left( e^{\hat{\mathcal{P}}(T-s)} e^{ij_n kX/\eta} - e^{ij_n\omega (T-s)/\varepsilon + ij_nkX/\eta  } 
               e^{\hat{\mathcal{Q}}_n(T-s)}\right) R_n(\hat{A}(s))ds \nonumber \\
&+& \sum_{j\notin J} H_j(T,X),
\label{4-20b}
\end{eqnarray}
where $H_j$ is defined by
\begin{eqnarray}
& & H_j(T,X) = \int^{T}_{0}\!  e^{\hat{\mathcal{P}}(T-s)} e^{ij\omega s/\varepsilon + ijkX/\eta}
                 C_j(\hat{A}(s))ds.
\label{4-14}
\end{eqnarray}

\textit{Step 3. estimate of $H_j$.}
Let $C([0,T_0];B^r)$ be a Banach space of functions $g(T,X)$ on $[0, T_0] \times \mathbf{R}^d$
such that $T \mapsto g(T, \, \cdot \,) \in B^r$ is continuous.
The norm is defined by
\begin{equation}
|| g ||_{C^{0,r}} := \max_{T\in [0,T_0]} || g(T, \,\cdot \,)  ||_{\eta}.
\label{4-15c}
\end{equation}
Let $C([0,T_0];B^r)^N$ be the product space with the norm defined by 
$\displaystyle || g ||_{C^{0,r}} = \max_{1\leq n\leq N}|| g_n ||_{C^{0,r}}$ for $g = (g_1, \cdots ,g_N)$.
Due to the existence theorem, a mild solution $\hat{A}(T,X)$ of (\ref{4-10b}) is included in $C([0,T_0];B^r)^N $.
The next lemma will be used several times.
\\[0.2cm]
\textbf{Lemma \thedef.} 
Suppose $f : B^r \to B^r$ is $C^1$.
There exists a function $h : [0, T_0] \times \mathbf{R}^d \times C([0,T_0];B^r)^N \to C([0,T_0];B^r)$,
which is bounded as $\eta \to 0$, such that
\begin{equation}
\sum_{j\notin J} H_j(T,X) = \eta h(T,X, \hat{A}(T,X)).
\label{4-16}
\end{equation}
Further, if $f : B^r \to B^r$ is $C^2$, $h(T,X, \hat{A})$ is Lipschitz continuous in $\hat{A} \in C([0,T_0];B^r)^N$.
\\[0.2cm]
\textbf{Proof.} 
In the $x$-coordinate, we have
\begin{eqnarray}
& & e^{\hat{\mathcal{P}}(T-s)}e^{ij\omega s/ \varepsilon + ijkX/\eta } C_j(\hat{A}(s,X)) \nonumber \\
 &=& \frac{1}{(2\pi)^d} \int\!  e^{ij\omega s/ \varepsilon + ijk(x+y)} C_j(\hat{A}(s,\eta x + \eta y))
\int \! e^{-iy\xi} e^{P(i\xi)(T-s)/\varepsilon }dy d\xi .
\label{4-17}
\end{eqnarray}
There exists a number $0<\theta <1$ such that
\begin{eqnarray*}
& & e^{\hat{\mathcal{P}}(T-s)}e^{ij\omega s/ \varepsilon + ijkX/\eta } C_j(\hat{A}(s,X)) \nonumber \\
 &=& \frac{1}{(2\pi)^d} \int\!  e^{ij\omega s/ \varepsilon + ijk(x+y)} C_j(\hat{A}(s,\eta x))
       \int \! e^{-iy\xi} e^{P(i\xi)(T-s)/\varepsilon }dy d\xi \\
&+& \frac{\eta }{(2\pi)^d} \int\!  e^{ij\omega s/ \varepsilon + ijk(x+y)}
    \sum^d_{i=1}\frac{\partial }{\partial y_i}\Bigl|_{y\mapsto \eta x+ \theta \eta y} C_j(\hat{A}(s,y)) \cdot y_i
      \int \! e^{-iy\xi} e^{P(i\xi)(T-s)/\varepsilon }dy d\xi.
\end{eqnarray*}
Since $\int \! e^{-iy\xi} e^{P(i\xi)(T-s)/\varepsilon }d\xi$ is rapidly decreasing in $y$,
the right hand side above exists.
We denote the first term and the second term above by $I_1$ and $I_2$, respectively;
$H_j = \int^{T}_{0}\! I_1ds + \int^{T}_{0}\! I_2ds$.
At first, we consider $I_2$.
Since $f$ is a $C^1$ function on $B^r$ and $[0,T_0]$ is a finite interval,
$f$ is regarded as a $C^1$ function on $C([0,T_0];B^r)$.
Since the derivatives $\partial C_j$'s are Fourier coefficients of $\partial f$, the series $\sum_{j\notin J} I_2$ converges and
there exists a function $h_2(T,X, \hat{A})$ such that
\begin{eqnarray*}
\sum_{j\notin J} \int^{T}_{0}\! I_2 ds = \eta h_2(T,X,\hat{A}).
\end{eqnarray*}
From the definition, we verify that $h_2$ is a mapping from $[0,T_0] \times \mathbf{R}^d \times C([0,T_0]; B^r)^N$
into $C([0,T_0]; B^0)$ (we will show later that this is a mapping into $C([0,T_0]; B^r)$).
Furthermore, if $f$ is $C^2$ so that $C_j$'s are $C^2$, then
$h_2(T,X, \hat{A})$ is $C^1$ in $\hat{A} \in C([0,T_0];B^r)^N$ (in particular, Lipschitz continuous).

Next, let us calculate the first term $I_1$.
Note that the equality
\begin{equation}
\int \!\int \! e^{ijky} e^{-iy\xi} e^{P(i\xi)(T-s)/\varepsilon }dy d\xi
 = (2\pi)^d e^{P(ijk)(T-s)/\varepsilon }
\label{4-18b}
\end{equation}
holds, which can be proved by the same way as Lemma 3.4.
Thus we obtain
\begin{eqnarray*}
\int^{T}_{0}\! I_1 ds
 = e^{ijkx}e^{P(ijk)T/\varepsilon } \int^{T}_{0}\!
    C_j(\hat{A}(s, \eta x)) e^{(ij\omega - P(ijk))s/\varepsilon } ds.
\end{eqnarray*}
If $\hat{A}$ is differentiable in $s$ (i.e. when the initial condition is included in the domain
of $\mathcal{Q}_n$), then integration by parts proves that the above quantity is of $O(\varepsilon )$.
When $\hat{A}$ is not differentiable, we need further analysis.

Let $J'$ be the set of integers $j$ such that $j\notin J$ and $\mathrm{Re}[P(ijk)] =0$.
Due to the assumption (D3), $J'$ is a finite set.
Put $ij\omega  - P(ijk) = p_j + iq_j$ with $p_j, q_j\in \mathbf{R}$.
When $j\notin J\cup J'$ (i.e. $p_j \neq 0$), 
the mean value theorem proves that there exists $0\leq \tau_j \leq T$ such that
\begin{eqnarray*}
\int^{T}_{0}\! C_j(\hat{A}(s, \eta x)) e^{(p_j + iq_j)s/\varepsilon } ds
&=& C_j(\hat{A}(\tau_j, \eta x)) e^{iq_j \tau _j / \varepsilon } \int^{T}_{0}\! e^{p_j s/ \varepsilon }ds \\
&=& \frac{\varepsilon }{p_j} C_j(\hat{A}(\tau_j, \eta x)) e^{iq_j \tau _j / \varepsilon }
    ( e^{p_j T / \varepsilon } - 1).
\end{eqnarray*}
Since $C_j$'s are Fourier coefficients of a $C^1$ function $f$ and 
since $p_{j} \to \infty$ as $|j| \to \infty$, the following series
\begin{eqnarray*}
& & \varepsilon \sum_{j\notin J\cup J'} e^{ijkx}e^{P(ijk)T/\varepsilon } \frac{1}{p_j} 
          C_j(\hat{A}(\tau_j, \eta x)) e^{iq_j \tau _j / \varepsilon } ( e^{p_j T / \varepsilon } - 1) \\
&=& \varepsilon \sum_{j\notin J\cup J'} e^{ijkx}e^{i\,\mathrm{Im} [P(ijk)] T/\varepsilon } \frac{1}{p_j} 
          C_j(\hat{A}(\tau_j, \eta x)) e^{iq_j \tau _j / \varepsilon } ( 1 - e^{-p_j T/ \varepsilon })
\end{eqnarray*}
converges, and there exists a function $h_1(T,X, \hat{A})$ from $[0,T_0] \times \mathbf{R}^d \times C([0,T_0]; B^r)^N$
into $C([0,T_0]; B^0)$, which is $C^1$ in $\hat{A}$, such that
\begin{eqnarray*}
\sum_{j\notin J\cup J'} \int^{T}_{0}\! I_1 ds = \varepsilon h_1(T,X,\hat{A}).
\end{eqnarray*}
To estimate the case $j\in J'$, we need the next lemma.
\\[0.2cm]
\textbf{Lemma \thedef.} 
For any constants $c, t_1 > 0$ that are independent of $\varepsilon $, 
a mild solution of (\ref{4-10b}) satisfies $|| \hat{A}_n(T + c\varepsilon ) -   \hat{A}_n(T) ||_{\eta} 
\sim O(\eta)$ for $\varepsilon t_1 \leq T \leq T_0$ and $n=1,\cdots ,N$.
\\[0.2cm]
\textbf{Proof.} 
Recall $x = X/\eta,\, t = T/\varepsilon $.
In the $(t,x)$-coordinates, linear semigroups satisfy
\begin{eqnarray*}
e^{\hat{Q}_n (T + c\varepsilon )}v_n - e^{\hat{Q}_n T}v_n
= \frac{1}{(2\pi)^d} \int\! v_n(\eta x + \eta y) \int\! e^{-iy \xi}  (e^{Q_n(i\xi) (t+c)} - e^{Q_n(i\xi)t})dyd\xi.  
\end{eqnarray*}
Since 
\begin{eqnarray*}
\frac{1}{(2\pi)^d} \int\! v_n(\eta x) \int\! e^{-iy \xi} e^{Q_n(i\xi) t}dyd\xi = v_n(\eta x)
\end{eqnarray*}
for any $t$, there exists $0<\theta <1$ such that
\begin{eqnarray*}
e^{\hat{Q}_n (T + c\varepsilon )}v_n - e^{\hat{Q}_n T}v_n
= \frac{\eta}{(2\pi)^d} \int\! \sum^d_{i=1} \frac{\partial v_n}{\partial x_i}(\eta x + \theta \eta y)y_i 
      \int\! e^{-iy \xi} e^{Q_n(i\xi)t} (e^{Q_n(i\xi) c} - 1)dyd\xi,
\end{eqnarray*}
which is of order $O(\eta)$ uniformly in $x$ and $t_1 \leq t$.
Next, the derivatives satisfy
\begin{eqnarray*}
\eta^{|\alpha |} \frac{\partial ^\alpha }{\partial X^\alpha } 
     \left( e^{\hat{Q}_n (T + c\varepsilon )}v_n - e^{\hat{Q}_n (T)}v_n \right) \!=\! 
\frac{1}{(2\pi)^d}\!\! \int\!\! v_n(\eta x +\eta y) \!\!
\int\! (i\xi)^\alpha e^{-iy\xi} e^{Q_n(i\xi)t}(e^{Q_n(i\xi) c}\!-\!1)dyd\xi.  
\end{eqnarray*}
By the same calculation as above, it turns out that 
\begin{eqnarray}
|| e^{\hat{Q}_n (T + c\varepsilon )}v_n - e^{\hat{Q}_n T}v_n ||_{\eta} \sim O(\eta)
\label{4-19b}
\end{eqnarray}
for $t_1 \leq t $; that is, for $\varepsilon t_1 \leq T$.
Then, Eq.(\ref{4-10b}) yields
\begin{eqnarray*}
& & \hat{A}_n(T + c\varepsilon ) -   \hat{A}_n(T) 
= e^{\hat{Q}_n (T + c\varepsilon )}v_n - e^{\hat{Q}_n T}v_n 
      +\int^{T+c\varepsilon }_{T}\!  e^{\hat{Q}_n (T + c\varepsilon -s)}R_n(\hat{A}(s))ds \\
& &   + \int^{T-\varepsilon t_1}_{0}\!  ( e^{\hat{Q}_n (T + c\varepsilon -s)}-e^{\hat{Q}_n (T -s)})R_n(\hat{A}(s))ds
+\int^{T}_{T-\varepsilon t_1}\!  ( e^{\hat{Q}_n (T + c\varepsilon -s)}-e^{\hat{Q}_n (T -s)})R_n(\hat{A}(s))ds.
\end{eqnarray*}
Using (\ref{4-19b}), we obtain the lemma. \hfill $\blacksquare$
\\

Suppose $j\in J'$, so that $ij\omega -P(ijk) = iq_j$;
\begin{eqnarray*}
\int^{T}_{0}\!I_1ds=e^{ijkx}e^{P(ijk)T/\varepsilon }\int^{T}_{0}\!C_j(\hat{A}(s,\eta x))e^{iq_js/\varepsilon }ds.  
\end{eqnarray*}
Since the set $J$ consists of all integers satisfying $P(ijk) = ij\omega $ (the assumptions (D3)), $q_j \neq 0$.
Put $I_{3,j} = \int^{T}_{0}\!C_j(\hat{A}(s,\eta x))e^{iq_js/\varepsilon }ds$.
Changing the variable $s \mapsto s + \varepsilon \pi / q_j$ yields
\begin{eqnarray*}
I_{3,j} = - \int^{T+\varepsilon \pi / q_j}_{\varepsilon \pi / q_j}\!
      C_j(\hat{A}(s- \varepsilon \pi / q_j,\eta x))e^{iq_js/\varepsilon }ds.
\end{eqnarray*}
Hence, we obtain
\begin{eqnarray*}
2I_{3,j} &=& \int^{\varepsilon \pi / q_j}_{0}\! C_j(\hat{A}(s,\eta x))e^{iq_js/\varepsilon }ds
         - \int^{T+\varepsilon \pi /q_j}_{T}\! C_j(\hat{A}(s- \varepsilon \pi /q_j,\eta x))e^{iq_js/\varepsilon }ds \\
& & \quad + \int^{T}_{\varepsilon \pi /q_j}\! \left( C_j(\hat{A}(s,\eta x))
 - C_j(\hat{A}(s- \varepsilon \pi /q_j,\eta x))\right) e^{iq_js/\varepsilon }ds \\
&=& \varepsilon \int^{\pi / q_j}_{0}\! C_j(\hat{A}(\varepsilon s,\eta x))e^{iq_js}ds
       - \varepsilon \int^{T/\varepsilon +\pi /q_j}_{T/\varepsilon }\! 
             C_j(\hat{A}(\varepsilon s-\varepsilon \pi /q_j,\eta x))e^{iq_js}ds \\
& & \quad + \, \eta \cdot \int^{T}_{\varepsilon \pi /q_j}\! \frac{C_j(\hat{A}(s,\eta x))
 - C_j(\hat{A}(s- \varepsilon \pi /q_j,\eta x))}{\eta} e^{iq_js/\varepsilon }ds.
\end{eqnarray*}
Lemma 4.5 shows that
\begin{eqnarray*}
\tilde{h}_{3,j}(\hat{A}) := \frac{C_j(\hat{A}(s,\eta x))
 - C_j(\hat{A}(s- \varepsilon \pi /q_j,\eta x))}{\eta }
\end{eqnarray*}
defines a function from $C([0,T_0];B^r)^N$ into $C([0,T_0];B^r)$, which is bounded as $\eta \to 0$.
Hence, there exists a function 
$h_{3,j}(T,X, \hat{A})$ from $[0,T_0] \times \mathbf{R}^d \times C([0,T_0]; B^r)^N$
into $C([0,T_0]; B^0)$ such that
\begin{eqnarray*}
\int^{T}_{0}\! I_1 ds = \eta h_{3,j}(T,X,\hat{A}).
\end{eqnarray*}
Further, if $f$ and $C_j$ are $C^2$, $h_{3,j}(T,X, \hat{A})$ is $C^1$ in $\hat{A}$.
Therefore, putting $h = h_2 + \eta ^{M-1}h_1 + \sum_{j\in J'} h_{3,j}$ proves Eq.(\ref{4-16})
satisfying $h(T,X, \hat{A}) \in C([0,T_0];B^0)$.

Let us estimate the derivative.
Eq.(\ref{4-17}) yields
\begin{eqnarray}
& & \eta ^{|\alpha |} \frac{\partial ^\alpha }{\partial X^\alpha }
\left( e^{\hat{\mathcal{P}}(T-s)}e^{ij\omega s/ \varepsilon + ijkX/\eta } C_j(\hat{A}(s,X)) \right) \nonumber \\
 &=& \frac{1}{(2\pi)^d} \int\!  e^{ij\omega s/ \varepsilon + ijk(x+y)} C_j(\hat{A}(s,\eta x + \eta y))
\int \! (i\xi) ^\alpha e^{-iy\xi} e^{P(i\xi)(T-s)/\varepsilon }dy d\xi .
\end{eqnarray}
Repeating the same argument, it turns out that $h(T,X, \hat{A})\in C([0,T_0];B^r)$ if $\hat{A}$ is in 
$C([0,T_0];B^r)^N$.
This completes the proof of the lemma. \hfill $\blacksquare$
\\

Because of this lemma, there exists a positive number $D_2$ such that
\begin{equation}
|| \sum_{j\notin J} H_j(T,X)||_{\eta} \leq \eta D_2
\label{4-21b}
\end{equation}
holds for $0\leq T\leq T_0$.

\textit{Step 4. estimate of $F(T)$.}
By using Eq.(\ref{4-12}) and (\ref{4-21b}),
we can show there exists a positive constant $D_3$ such that $|| F(T) ||_\eta\leq \eta D_3$
for $\varepsilon t_0 \leq T\leq T_0$.
Therefore, we obtain
\begin{eqnarray*}
& & || \hat{u} - \sum^N_{n=1} \hat{A}_n e^{ij_n\omega T/\varepsilon + ij_nkX/\eta }  ||_\eta \\
& \leq & \eta D_3 + D_1L \int^{T}_{\varepsilon t_0}\! e^{D_1L(T-s)} \eta D_3 ds 
 + D_1L \int^{\varepsilon t_0}_{0}\! e^{D_1L(T-s)} || F(s) ||_{\eta} ds \sim O(\eta).
\end{eqnarray*}
for $\varepsilon t_0\leq T\leq T_0$.
Changing to the $(t,x)$-coordinate proves Theorem 4.2. \hfill $\blacksquare$
\\[0.2cm]
\textbf{Proof of Thm.4.3.}
Let us consider the systems (4.15).
Recall that in this situation, $C_j(A_1,\cdots ,A_N, A_{-1}, \cdots ,A_{-N})$ is defined by Eq.(\ref{4-35b}).
$S_n(A)$ is defined by
\begin{eqnarray*}
S_n(A_1 ,\cdots ,A_N) = R_n(A_1 ,\cdots ,A_N, A_1,\cdots ,A_N) = 
C_{j_n}(A_1 ,\cdots ,A_N, A_1,\cdots ,A_N) 
\end{eqnarray*} 
for $n=1,\cdots , N$.
Again we assume $D = d$ and use the same notation as the previous proof.
A mild solution of (\ref{4-39a}) written in the $(T,X)$-coordinate satisfies
\begin{equation}
\hat{u} = e^{\hat{\mathcal{P}}T} (\hat{u}(0))
        + \int^{T}_{0}\! e^{\hat{\mathcal{P}}(T-s)} f(\hat{u}(s))ds,
\label{4-29b}
\end{equation}
with the initial condition $\hat{u}(0,X) = \hat{u}(0)$.
Let us consider the system of $2N$-integral equations of $w_+ := (w_1, \cdots ,w_N)$ 
and $w_- := (w_{-1}, \cdots ,w_{-N})$ of the form
\begin{eqnarray}
e^{ij_n \omega T/\varepsilon + ij_n kX/\eta} w_n 
&=& e^{\hat{\mathcal{P}}T} (e^{ij_n kX/\eta}w_n(0))
 +\int^{T}_{0}\! e^{\hat{\mathcal{P}}(T-s)}R_n(w_+(s),w_-(s))e^{ij_n \omega s/\varepsilon +ij_n kX/\eta}ds \nonumber \\
& & +\frac{1}{2N}\sum_{j\notin  J} 
      \int^{T}_{0}\! e^{\hat{\mathcal{P}}(T-s)} C_j(w_+(s),w_-(s))e^{ij \omega s/\varepsilon + ij kX/\eta} ds, \\
e^{-ij_n \omega T/\varepsilon - ij_n kX/\eta} w_{-n} 
&=& e^{\hat{\mathcal{P}}T} (e^{-ij_n kX/\eta}w_{-n}(0))
 +\int^{T}_{0}\! e^{\hat{\mathcal{P}}(T-s)}R_{-n}(w_+(s),w_-(s))e^{-ij_n \omega s/\varepsilon -ij_n kX/\eta}ds \nonumber \\
& & +\frac{1}{2N}\sum_{j\notin  J} 
      \int^{T}_{0}\! e^{\hat{\mathcal{P}}(T-s)} C_{-j}(w_+(s),w_-(s))e^{-ij \omega s/\varepsilon - ij kX/\eta} ds,
\end{eqnarray}
which has a unique solution $w_n \in C([0,T_0];B^r)$ satisfying $w_n(0,X) = w_n(0) \in B^r$.
This yields
\begin{eqnarray*}
& & \sum^N_{n=-N}e^{ij_n \omega T/\varepsilon + ij_n kX/\eta} w_n \\
&=& e^{\hat{\mathcal{P}}T} (\sum^N_{n=-N}e^{ij_n kX/\eta} w_n(0))
    + \sum^\infty_{j=-\infty} \int^{T}_{0}\!
          e^{\hat{\mathcal{P}}(T-s)} C_j(w_+(s),w_-(s))e^{ij \omega s/\varepsilon + ij kX/\eta}  ds \\
&=& e^{\hat{\mathcal{P}}T} (\sum^N_{n=-N}e^{ij_n kX/\eta}w_n(0))
    + \int^{T}_{0}\! e^{\hat{\mathcal{P}}(T-s)}
         f(\textstyle \sum^N_{n=-N} w_n(s)e^{ij_n \omega s/\varepsilon + ij_n kX/\eta}  )ds,
\end{eqnarray*}
where we used the abbreviation $\sum^N_{n=-N} := \sum^{-N}_{n=-1} + \sum^N_{n=1}$.
Recall also that $j_{-n} = -j_n$.

\noindent This means that $\hat{u} = \sum^N_{n=-N}e^{ij_n \omega T/\varepsilon + ij_n kX/\eta} w_n$ 
when $\hat{u}(0) = \sum^N_{n=-N}e^{ij_n kX/\eta} w_n(0)$.

By using (D4), we can show the equalities
\begin{eqnarray*}
\overline{e^{\hat{\mathcal{P}}T}(e^{ijkX/\eta} w_n)} = e^{\hat{\mathcal{P}}T}(e^{-ijkX/\eta} \overline{w}_n),
\end{eqnarray*} 
and
\begin{eqnarray*}
\overline{C_j(A_1, \cdots ,A_N, A_{-1}, \cdots ,A_{-N})} = 
C_{-j}(\overline{A}_{-1}, \cdots , \overline{A}_{-N}, \overline{A}_1, \cdots , \overline{A}_{N}).
\end{eqnarray*} 
Therefore, $\{  w_n = w_{-n} \in \mathbf{R} \}_{n=1}^N$ is the invariant set of the $2N$-equations.
By putting $w_n = w_{-n} \in \mathbf{R}$ and $w:= w_+ = w_- $,
the system is reduced to $N$-equations of the form
\begin{eqnarray}
e^{ij_n \omega T/\varepsilon + ij_n kX/\eta} w_n 
&=& e^{\hat{\mathcal{P}}T} (e^{ij_n kX/\eta}w_n(0))
 +\int^{T}_{0}\! e^{\hat{\mathcal{P}}(T-s)}S_n(w(s))e^{ij_n \omega s/\varepsilon +ij_n kX/\eta}ds \nonumber \\
& & +\frac{1}{2N}\sum_{j\notin  J} 
      \int^{T}_{0}\! e^{\hat{\mathcal{P}}(T-s)} C_j(w(s))e^{ij \omega s/\varepsilon + ij kX/\eta} ds, 
\label{4-15}
\end{eqnarray}
for $n=1,\cdots ,N$, where $ C_j(w):= C_j(w,w)$.

Suppose that there exists $\phi = (\phi_1, \cdots ,\phi_N) \in \mathbf{R}^N$ satisfying
the assumptions of Thm.4.3.
Without loss of generality, we assume that the Jacobi matrix of $(S_1, \cdots ,S_N)$ at $A = \phi$ is diagonal.
Thus we put
\begin{equation}
S_n(\phi_1, \cdots ,\phi_N) = 0, \quad \frac{\partial S_n}{\partial x_m}(\phi_1, \cdots ,\phi_N)
 = -\beta_n \cdot \delta _{n,m},\quad \mathrm{Re}[\beta_n] > 0,
\end{equation}
for $n,m=1,\cdots ,N$.
Put $w_n = \phi_n + \eta W_n$.
Due to the assumption of Thm.4.3, $S_n$ is $C^2$ having the Lipschitz continuous second derivatives.
Hence, there is a Lipschitz continuous function $\hat{S}_n(\, \cdot \,, \eta) : C([0,T_0];B^r)^N \to C([0,T_0];B^r)$,
which is bounded as $\eta \to 0$, such that
\begin{equation}
S_n(\phi + \eta W) = -\eta \beta_n W_n + \eta^2 \hat{S}_n(W, \eta),
\quad W = (W_1, \cdots ,W_n).
\end{equation}
We denote $\hat{S}_n(W, \eta)$ by $\hat{S}_n(W)$ for simplicity.
Note that Eq.(\ref{4-18b}) gives
\begin{eqnarray*}
e^{\hat{\mathcal{P}}T}(e^{ij_nkX/\eta} \phi_n)
 = e^{P(ij_nk)T/\varepsilon } e^{ij_nkX/\eta} \phi_n = e^{ij_n \omega T/\varepsilon + ij_n kX/\eta}\phi_n.
\end{eqnarray*}
Thus Eq.(\ref{4-15}) is rewritten as 
\begin{eqnarray*}
\!\!& & e^{ij_n \omega T/\varepsilon + ij_n kX/\eta} W_n 
= e^{\hat{\mathcal{P}}T} (e^{ij_n kX/\eta}W_n(0))
    -\beta_n \int^{T}_{0}\! e^{\hat{\mathcal{P}}(T-s)}W_n(s) e^{ij_n \omega s/\varepsilon +ij_nkX/\eta } ds \\
\!\!& +&\!\!\!\! \eta \!\!
\int^{T}_{0}\!\! e^{\hat{\mathcal{P}}(T-s)}\hat{S}_n(W(s)) e^{ij_n \omega s/\varepsilon + ij_nkX/\eta } ds 
 + \frac{1}{2 \eta N}\sum_{j\notin  J} \!\int^{T}_{0}\! \!
         e^{\hat{\mathcal{P}}(T-s)} C_j(\phi + \eta W(s)) e^{ij \omega s/\varepsilon +ijkX/\eta } ds.
\end{eqnarray*}
Remark that the second term in the right hand side is linear in $W_n$.
Therefore, we can show that this equation is rewritten as
\begin{eqnarray}
 e^{ij_n \omega T/\varepsilon + ij_n kX/\eta} W_n 
&=& e^{(\hat{\mathcal{P}} - \beta_n)T} (e^{ij_n kX/\eta}W_n(0))
 + \eta \int^{T}_{0}\! e^{(\hat{\mathcal{P}} - \beta_n)(T-s)}
         \hat{S}_n(W(s)) e^{ ij_n \omega s/\varepsilon +ij_nkX/\eta } ds \nonumber \\
& &  + \frac{1}{2\eta N}\sum_{j\notin  J} \int^{T}_{0}\! 
         e^{(\hat{\mathcal{P}} - \beta_n)(T-s)} C_j(\phi + \eta W(s)) e^{ij \omega s/\varepsilon +ijkX/\eta } ds.
\label{4-18}
\end{eqnarray}
Motivated by this equation, let us consider the system of integral equations of the form
\begin{eqnarray}
 e^{ij_n \omega T/\varepsilon + ij_n kX/\eta} W_n^*(T) 
&=& \eta \int^{T}_{-\infty}\!  
      e^{(\hat{\mathcal{P}} - \beta_n)(T-s)} e^{ij_n \omega s/\varepsilon + ij_nkX/\eta } \hat{S}_n(W^*(s)) ds \nonumber \\
& & +   \frac{1}{2\eta N}\sum_{j\notin  J} \int^{T}_{-\infty}\! 
    e^{(\hat{\mathcal{P}} - \beta_n)(T-s)} e^{ij \omega s/\varepsilon + ijkX/\eta } C_j(\phi + \eta W^*(s)) ds. \quad\quad
\label{4-19}
\end{eqnarray}
We will show later that a solution $W^* = (W_1^*,\cdots ,W_N^*)$ of this system satisfies Eq.(\ref{4-18}) with a suitable 
initial condition $W_n(0) = W_n(0,X)$.

To prove the existence of a periodic solution, let $C^r_p$ be the set of functions $\psi (T,X)$ in 
$C([0,T_0]; B^r)$ such that
\begin{equation}
 \left\{ \begin{array}{l}
\displaystyle 
\psi (T + \frac{2\pi \varepsilon }{\omega }, X_1, \cdots , X_d) 
 = \psi (T, X_1, \cdots , X_d), \quad \text{when $\omega \neq 0$},  \\
\psi (T, X_1, \cdots , X_d)\,\,\, \text{is constant in $T$ when $\omega =0$},  \\
\end{array} \right.
\label{4-20}
\end{equation}
and
\begin{equation}
\left\{ \begin{array}{l}
\displaystyle 
\psi (T , X_1, \cdots ,X_{j-1},X_j + \frac{2\pi \eta}{k_j}, X_{j+1},\cdots , X_d) 
 = \psi (T, X_1, \cdots , X_d), \quad \text{when $k_j \neq 0$},  \\
\psi (T, X_1, \cdots , X_d)\,\,\, \text{is constant in $X_j$ when $k_j =0$},  \\
\end{array} \right.
\label{4-21}
\end{equation}
for $j= 1,\cdots ,d$.
By the norm $|| \psi ||_{C^r_p}:= \max_{T\in \mathbf{R}} || \psi (T, \cdot \,) ||_\eta$, 
$C^r_p$ becomes a Banach space, which is a closed subspace of $C([0,T_0]; B^r)$.
Let $(C^r_p)^N$ be the product space
with the norm $|| \psi ||_{C^r_p}= \max_{1\leq n\leq N}||\psi_n||_{C^r_p}$ for $\psi = (\psi_1, \cdots ,\psi_N)$.
Define mappings $\Omega _{1,n}$ and $\Omega _{2,n}$ to be
\begin{eqnarray*}
 (\Omega _{1,n}W)(T,X)
& =& \eta \int^{T}_{-\infty}\!  e^{-ij_n \omega (T-s) /\varepsilon - ij_n kX/\eta} 
         e^{(\hat{\mathcal{P}} - \beta_n)(T-s)} e^{ij_nkX/\eta } \hat{S}_n(W(s,X)) ds, \nonumber \\
(\Omega _{2,n}W)(T,X) &=& \!\! \frac{1}{2\eta N}\! \sum_{j\notin  J} \!\int^{T}_{-\infty}\!\! 
        e^{-ij_n\omega T/\varepsilon - ij_nkX/\eta }  e^{(\hat{\mathcal{P}} - \beta_n)(T-s)} 
        e^{ij\omega s/\varepsilon + ijkX/\eta } C_j(\phi \!+\! \eta W(s,X)) ds,
\end{eqnarray*}
for $n= 1,\cdots ,N$.
\\[0.2cm]
\textbf{Lemma \thedef.}
$\Omega _{1,n}$ and $\Omega _{2,n}$ are mappings from $(C^r_p)^N$ into $C^r_p$.
\\[0.2cm]
\textbf{Proof.}
By using the expression (\ref{3-5b}) of the semigroup, we can show that if $W\in (C^r_p)^N$, 
then there exists a positive constant $D_1$ such that
\begin{eqnarray*}
\eta^{|\alpha |} \Bigl| \frac{\partial^\alpha }{\partial X^\alpha }
   e^{-ij_nkX/\eta }e^{\hat{\mathcal{P}}(T-s)}e^{ij_nkX/\eta } \hat{S}_n(W(s,X)) \Bigr| \leq D_1
\end{eqnarray*}
for $ X\in \mathbf{R}^d,\, -\infty< s \leq T$ and $|\alpha |=1,\cdots ,r$.
Thus we obtain
\begin{eqnarray*}
|| \Omega _{1,n}W||_{C^r_p}
\leq \eta \cdot \sup_{T\in \mathbf{R}}\int^{T}_{-\infty}\!  D_1 e^{-\beta_n(T-s)} ds.
\end{eqnarray*}
Since $\mathrm{Re}[\beta_n] > 0$, the right hand side above exists.
The periodicity conditions (\ref{4-20}),(\ref{4-21}) immediately follow from the definition.
The proof for $\Omega _{2,n}$ is done in the same way. \hfill $\blacksquare$
\\[0.2cm]
\textbf{Lemma \thedef.}
Fix a positive number $\delta $ and let $\mathcal{D} = \{ W\in (C^r_p)^N \, | \, 
|| W ||_{C^r_p} \leq \delta  \}$ be a closed ball in $(C^r_p)^N$.
If $\eta = \eta (\delta ) > 0$ is sufficiently small, 
$(\Omega _{1,1}, \cdots ,\Omega _{1,N})$ and $(\Omega _{2,1}, \cdots ,\Omega _{2,N})$
are contraction mappings on $\mathcal{D}$. 
\\[0.2cm]
\textbf{Proof.}
It follows from the proof of Lemma 4.6 that $|| \Omega _{1,n}W||_{C^r_p}$ is of order $O(\eta )$.
Thus it is easy to verify that $(\Omega _{1,1}, \cdots ,\Omega _{1,N})$ is a contraction mapping on $\mathcal{D}$
if $\eta$ is sufficiently small.
Next, by the same way as the proof of Lemma 4.4, we can show that there exists a function
$h_n : \mathbf{R} \times \mathbf{R}^d \times (C^r_p)^N \to C^r_p$ such that
\begin{equation}
(\Omega _{2,n}W)(T, X) = \frac{1}{2\eta N} \cdot \eta h_n(T,X, \phi + \eta W (T,X)),
\end{equation}
where $h_n(T, X, \,\cdot \,)$ is Lipschitz continuous.
Therefore, there exists $L_n > 0$ such that
\begin{eqnarray*}
|| \Omega _{2,n}W - \Omega _{2,n}V ||_{C^r_p}
\leq \frac{L_n}{2N} || \eta W - \eta V ||_{C^r_p} \sim O(\eta ).
\end{eqnarray*}
This proves that $(\Omega _{2,1}, \cdots ,\Omega _{2,N})$ is contraction if $\eta$ is sufficiently small.
\hfill $\blacksquare$
\\

Due to this lemma, the system (\ref{4-19}) of integral equations has a unique solution
$W^* = (W^*_1, \cdots ,W^*_N)$ in $(C^r_p)^N$, which is periodic in $X$ and $T$.
\\[0.2cm]
\textbf{Lemma \thedef.} 
The solution $W^*$ is a solution of (\ref{4-18}) satisfying the initial condition
\begin{eqnarray}
e^{ij_nkX /\eta}W_n(0) 
&=& \eta \int^{0}_{-\infty}\! e^{-(\hat{\mathcal{P}} - \beta_n)s} e^{ij_n \omega s/ \varepsilon + ij_nkX/\eta } 
       \hat{S}_n(W^*(s)) ds \nonumber \\
&+ &   \frac{1}{2\eta N}\sum_{j\notin  J} \int^{0}_{-\infty}\! 
     e^{-(\hat{\mathcal{P}} - \beta_n)s} e^{ij \omega s/ \varepsilon +ijkX/\eta } C_j(\phi + \eta W^*(s)) ds.
\label{4-22}
\end{eqnarray}
\textbf{Proof.}
This follows from the substitution of (\ref{4-22}) into (\ref{4-18}). \hfill $\blacksquare$
\\

Now we have proved that the system (\ref{4-15}) has a solution $w_n(T,X) = \phi_n + \eta W_n^*(T,X)$
satisfying $W_n^* \in C^r_p$.
Therefore, the equation (\ref{4-29b}) has a solution
\begin{equation}
\hat{u}(T,X) = \sum^N_{n=1} \left( \phi_n + \eta W_n^*(T, X) \right) \cdot 
(e^{ij_n\omega T/\varepsilon + ij_nkX/\eta} + e^{-ij_n\omega T/\varepsilon - ij_nkX/\eta}).
\end{equation}
Changing to the $(t,x)$-coordinate yields a mild solution of (\ref{4-39a}) of the form
\begin{equation}
u(t,x) = \sum^N_{n=1} \left( \phi_n + \eta W_n^*(\varepsilon t, \eta x) \right) \cdot
(e^{ij_n\omega t + ij_nkx} + e^{-ij_n\omega t - ij_nkx}),
\end{equation}
which proves the first part of Thm.4.3.

Finally, let us prove the stability part of Thm.4.3.
Let $W^*(T, X)$ be the periodic solution of (\ref{4-18}).
There exists a positive constant $D_1\geq 1$ such that 
$||e^{-ij_nkX/\eta} e^{\hat{\mathcal{P}}T} e^{ij_nkX/\eta}||_\eta \leq D_1$ for any $T>0$ and $n=1,\cdots ,N$.
Fix a positive number $M$ and put $\delta = M/(2D_1)$.
Let $\mathcal{D} = \{ W\in (B^r)^N \, | \, || W-W^*(0, \cdot \,) ||_{\eta} \leq \delta \}$
be a neighborhood of the periodic solution $W^*(0, X)$ at $T = 0$.
Due to the existence theorem of mild solutions, there exists $T_0 > 0$ such that when $W(0) \in \mathcal{D}$,
then Eq.(\ref{4-18}) has a solution $W(T)$ in $(B^r)^N$ for $0\leq T\leq T_0$.
We define a time $T_0$ map as
\begin{equation}
\mathcal{T} : \mathcal{D} \to (B^r)^N, \quad W(0)\mapsto W(T_0).
\end{equation}
\textbf{Lemma \thedef.} If $\eta > 0$ is sufficiently small, $T_0$ can be taken so that
$\mathcal{T}$ is a mapping on $\mathcal{D}$.
\\[0.2cm]
\textbf{Proof.}
Define $\Omega _{3,n}$ and $\Omega _{4,n}$ to be
\begin{eqnarray*}
 (\Omega _{3,n}W)(T,X)
& =& \int^{T}_{0}\!  e^{-ij_n \omega (T-s) /\varepsilon - ij_n kX/\eta} 
         e^{(\hat{\mathcal{P}} - \beta_n)(T-s)} e^{ij_nkX/\eta } \hat{S}_n(W(s,X)) ds, \nonumber \\
(\Omega _{4,n}W)(T,X) &=& \!\! \frac{1}{2 N}\! \sum_{j\notin  J} \!\int^{T}_{0}\!\! 
        e^{-ij_n\omega T/\varepsilon - ij_nkX/\eta }  e^{(\hat{\mathcal{P}} - \beta_n)(T-s)} 
        e^{ij\omega s/\varepsilon + ijkX/\eta } C_j(\phi \!+\! \eta W(s,X)) ds,
\end{eqnarray*}
for $n= 1,\cdots ,N$.
Eq.(\ref{4-18}) gives
\begin{eqnarray*}
W_n - W^*_n &=& e^{-ij_n \omega T/ \varepsilon -ij_nkX/\eta}
                  e^{(\hat{\mathcal{P}} - \beta_n)T}e^{ij_nkX/\eta} (W_n(0) - W_n^*(0)) \nonumber \\
& & \quad \quad + \eta (\Omega _{3,n}W - \Omega _{3,n}W^*) + \frac{1}{\eta}(\Omega _{4,n}W - \Omega _{4,n}W^*).
\end{eqnarray*}
By the same way as the proof of Lemma 4.4, we can prove that there exist Lipschitz continuous functions
$h_n : [0,T_0] \times \mathbf{R}^d \times C([0,T_0],B^r)^N \to C([0,T_0],B^r)$ such that 
$\Omega _{4,n}W = \eta h_n(T,X, \phi + \eta W)$.
This provides
\begin{eqnarray*}
W_n - W^*_n &=&
  e^{-ij_n \omega T/ \varepsilon -ij_nkX/\eta}e^{(\hat{\mathcal{P}} - \beta_n ) T} e^{ij_nkX/\eta}(W_n(0) - W_n^*(0)) \nonumber \\
& & \quad \quad + \eta (\Omega _{3,n}W - \Omega _{3,n}W^*) +  h_n(T,X, \phi + \eta W)- h_n(T,X, \phi + \eta W^*).
\end{eqnarray*}
Since $|| W(0) - W^*(0) ||_\eta \leq \delta = M/(2D_1)$ for $W(0) \in \mathcal{D}$, 
we can assume that $T_0$ is chosen so that $|| W(T) - W^* (T) ||_\eta \leq M$ for $0\leq T\leq T_0$.
Put $\beta := \min_{1\leq n\leq N} \beta_n$.
Since $\hat{S}_n$ and $h_n$ are locally Lipschitz continuous, there exist positive constants $D_2,D_3$ such that
\begin{eqnarray*}
 e^{\beta T} || W(T)-W^*(T) ||_\eta 
\leq M/2 + \eta D_2 \int^{T}_{0}\! e^{\beta s} || W(s)-W^*(s) ||_\eta ds
 + \eta D_3 e^{\beta T} || W-W^* ||_{C^{0,r}}
\end{eqnarray*}
for $0\leq T\leq T_0$ (see Eq.(\ref{4-15c}) for the definition of the norm $|| \cdot ||_{C^{0,r}}$ on
$C([0,T_0]; B^r)^N$).
The Gronwall inequality gives
\begin{eqnarray*}
& &  e^{\beta T} || W(T)-W^* (T)||_\eta \\
& \leq & M/2  + \eta D_3 e^{\beta T} || W-W^* ||_{C^{0,r}} 
+\eta D_2 \int^{T}_{0}\! e^{\eta D_2 (T-s)} \left(  M/2 + \eta D_3 e^{\beta s} || W-W^* ||_{C^{0,r}} \right) ds \\
&=& M e^{\eta D_2T}/2 + \frac{\eta D_3 \beta}{\beta - \eta D_2}e^{\beta T} || W-W^* ||_{C^{0,r}} 
 - \frac{\eta ^2 D_2D_3}{\beta - \eta D_2} e^{\eta D_2 T} || W-W^* ||_{C^{0,r}} 
\end{eqnarray*}
for $0\leq T\leq T_0$.
Hence, we obtain
\begin{eqnarray*}
& &  || W(T)-W^*(T) ||_\eta 
\leq  M e^{(\eta D_2 - \beta )T}/2
       + \frac{\eta D_3 \beta}{\beta - \eta D_2} || W-W^* ||_{C^{0,r}}.
\end{eqnarray*}
Using the standard existence theorem, we can verify that
a solution $W_n(T,X)$ of (\ref{4-18}) is bounded as $\eta \to 0$.
Thus, $|| W-W^* ||_{C^{0,r}}$ is bounded as $\eta \to 0$.
Therefore, if $\eta$ is sufficiently small, we obtain $|| W(T)-W^*(T) ||_\eta \leq M/2$.
This implies that $T_0$ can be taken arbitrarily large and $|| W(T)-W^*(T) ||_\eta \leq M$ holds for any $T>0$.
Hence, if $T_0$ is sufficiently large, we obtain $|| W(T_0)-W^*(T_0) ||_\eta \sim O(\eta)$.
Since $W^*$ is $2\pi \varepsilon /\omega $-periodic in $T$, we can choose $T_0$
so that $|| W(T_0)-W^*(0) ||_\eta \sim O(\eta)$, which proves $W(T_0) \in \mathcal{D}$. \hfill $\blacksquare$
\\[0.2cm]
\textbf{Lemma \thedef.} If $\eta > 0$ is sufficiently small and $T_0$ is sufficiently large,
$\mathcal{T}$ is a contraction mapping on $\mathcal{D}$.
\\[0.2cm]
\textbf{Proof.}
Let $W(T)$ and $V(T)$ be two solutions of (\ref{4-18}) with the initial conditions
$W(0), V(0) \in \mathcal{D}$, respectively.
By the same calculation as above, we can show the inequality
\begin{eqnarray*}
 || W(T)-V(T) ||_\eta 
\leq  D_1 || W(0)-V(0) ||_{\eta} e^{(\eta D_2 - \beta )T} + \eta D_4 || W-V||_{C^{0,r}},
\end{eqnarray*}
where we put $D_4 = D_3\beta /(\beta - \eta D_2)$.
Since $|| W-V ||_{C^{0,r}} = \max_{0\leq T \leq T_0}|| W(T)-V(T) ||_\eta$, we obtain
\begin{eqnarray*}
 || W-V ||_{C^{0,r}} \leq  D_1 || W(0)-V(0) ||_{\eta} + \eta D_4 || W-V||_{C^{0,r}} . 
\end{eqnarray*}
Substituting this into the above inequality provides
\begin{eqnarray*}
|| W(T)-V(T) ||_\eta \leq 
 \left( D_1e^{(\eta D_2 - \beta )T} + \frac{\eta D_1D_4}{1- \eta D_4} \right) || W(0)-V(0) ||_{\eta},
\end{eqnarray*}
which proves the lemma. \hfill $\blacksquare$
\\

Take $T_0 = 2\pi \varepsilon /\omega \cdot l$ so that Lemmas 4.9 and 4.10 hold,
where $l$ is a sufficiently large integer.
Due to the lemma, there exists a unique function $W^{**}(X) \in \mathcal{D}$ such that any solutions
$W(T,X)$ in $\mathcal{D}$ satisfy $W(nT_0, \,\cdot \,) \to W^{**}(\cdot)$ as $n\to \infty$.
If $W^{**} \neq W^*$, then $W^*(nT_0, \,\cdot \,)$ converges to $W^{**}$ as $n\to \infty$.
This contradicts with the fact that $W^*(nT_0, \,\cdot \,) = W^* (0, \, \cdot \,)$ is independent of $n$.
Hence, $W^{**} = W^*$ and $W(T,X)$ in $\mathcal{D}$ converges to $W^*$ as $T\to \infty$.
Since a mild solution of (\ref{4-39a}) is written as 
\begin{equation}
u(t,x) = \sum^N_{n=1} \left( \phi_n + \eta W_n(\varepsilon t, \eta x) \right) 
\cdot (e^{ij_n\omega t + ij_nkx}+e^{-ij_n\omega t - ij_nkx}),
\end{equation}
the proof of Thm.4.3 is completed. \hfill $\blacksquare$

%%%%%%%%%%%%%%%%%%%%%%%%%%%%%%%%%%%%%%%%%%%%%%%%%%%%%%%%%%%%%%%%%%%%%%%%%%%%%%%%%%%%%%%%%%%%%%%%%%

\subsection{Higher dimensional case}

Suppose $u = (u_1, \cdots ,u_m) \in \mathbf{C}^m$ and $x = (x_1, \cdots ,x_d) \in \mathbf{R}^d$.
For fixed $1\leq D\leq d$, we use the same notation $x = (\hat{x}_1, \hat{x}_2),\,
\alpha = (\beta, \gamma )$ as in Sec.4.1.
Let $\{ P_{ij}(x) \}_{i,j=1}^m$ be the set of polynomials of $x$.
The $m\times m$ matrix $P(x)$, the differential operator $\mathcal{P}$ and 
eigenvalues $\lambda _1(\xi) , \cdots ,\lambda _m(\xi)$ are defined in the same way as Sec.3.2.
We suppose for simplicity that only $\lambda _1(\xi)$ contributes to the center subspace of $\mathcal{P}$
(see (E1) below).
Extending to more general situations is not difficult (see Remark 3.7 and Example 4.12).
\\[0.2cm]
\textbf{(E0)} The matrix $P(i\xi)$ is diagonalizable for any $\xi \in \mathbf{R}^d$.
\\
\textbf{(E1)} $\mathrm{Re}[\lambda _1(\xi)] \leq 0$ and $\mathrm{Re}[\lambda _j(\xi)] < 0$
for any $\xi \in \mathbf{R}^d$ and $j=2,\cdots , m$.
\\
\textbf{(E2)}\, There exist $\omega \in \mathbf{R},\, k\in \mathbf{R}^d\, ((\omega ,k) \neq (0,0))$, a finite set of integers 
$J = \{ j_1, \cdots ,j_N\}$ and $\{ M_1, \cdots ,M_N\}$ such that
\begin{eqnarray*}
& & \lambda _1(j_nk) = ij_n\omega, \,\, (n=1,\cdots ,N),\\
& & \frac{\partial ^\beta \lambda_1 }{\partial \hat{x}_1^\beta }(j_nk) = 0, \,\,\text{for any $\beta $ such that}
\,\, |\beta | = 1, \cdots  ,M_n-1, \,\, (n=1,\cdots ,N) ,\\
& & \frac{\partial ^\beta \lambda _1}{\partial \hat{x}_1^\beta }(j_nk) \neq 0,\,\, \text{for some $\beta _n$ such that}
\,\, |\beta _n| = M_n,\,\, (n=1,\cdots ,N).
\end{eqnarray*}
The set $J$ consists of \textit{all} integers satisfying $\lambda _1(jk) = ij\omega $.
\\
\textbf{(E3)} For $n=1, \cdots ,N$, define $Q_n(x)$ and $\mathcal{Q}_n$ by
\begin{equation}
Q_n(x) = Q_n(\hat{x}_1, 0) = \sum_{|\beta | = M_n} \frac{1}{(\beta _1!)\cdots (\beta _D !)} 
\frac{\partial ^\beta \lambda _1}{\partial \hat{x}_1^\beta }(j_nk) (\hat{x}_1 / i)^\beta , \quad 
\mathcal{Q}_n = Q_n(\partial _1 , \cdots ,\partial _D, 0 ,\cdots ,0).
\end{equation}
Then, both of $\mathcal{P}$ and $\mathcal{Q}_n$ are elliptic in the sense that
there exist $c_1,c_2 > 0$ such that $\mathrm{Re}[\lambda _j(\xi )] < -c_2|\xi|^2 \,( j=1,\cdots ,m)$  and 
$\mathrm{Re}[Q_n(i\hat{\xi}_1)] < -c_2|\hat{\xi}_1|^2 \,( n=1,\cdots ,N)$ hold for $|\xi|, |\hat{\xi}_1| \geq c_1$,
where $\hat{\xi}_1 = (\xi_1, \cdots ,\xi_D)$.
\\

When $m=1$, $\lambda _1(\xi) = P(i\xi)$, so that the above assumptions and $Q_n$ are reduced to those given in Sec.4.1.
Let $B^r = BC^r(\mathbf{R}^r; \mathbf{C})$ and $(B^r)^m = B^r \times \cdots  \times B^r$ a product space.
The norm on $(B^r)^m$ is defined by $|| u || = \max_{1\leq n\leq m}|| u_n ||$.
Note that $\mathcal{P}$ is an operator densely defined on $(B^r)^m$, while $\mathcal{Q}_n$ is an operator
densely defined on $B^r$.

Let $\bm{w_n} = (w_{n,1}, \cdots ,w_{n,m})$ be an eigenvector of $P(ij_nk)$ associated with 
the eigenvalue $\lambda _1(j_nk)= ij_n \omega $ for $n=1, \cdots ,N$.
The projection to the eigenspace $\mathrm{span} \{ \bm{w}_n\}$ is denoted by $\Pi_n$,
and the projection to the eigenspace associated with the other eigenvalues 
$\lambda _2(j_nk), \cdots , \lambda _m(j_nk)$ is denoted by $\Pi_n^\bot = id - \Pi_n$.
Functions $e^{ij_1kx}\bm{w}_1,\cdots ,e^{ij_Nkx}\bm{w}_N$ span the center subspace of $\mathcal{P}$.

When $k=(k_1, \cdots ,k_d) \neq 0$, we can assume without loss of generality that $k_1\neq 0$ and $j_1 = 1$.
For a given function $f : (B^r)^m \to (B^r)^m$, define a function $\bm{C}_j : (B^r)^N \to (B^r)^m$ by
\begin{equation}
\bm{C}_j (A) = \bm{C}_j(A_1 , \cdots ,A_N)
 = \frac{k_1}{2\pi}\int^{2\pi / k_1}_{0}\! f(\sum^N_{n=1} A_n e^{ij_n k_1x_1}\bm{w}_n) e^{-ijk_1x_1} dx_1. 
\label{4-47}
\end{equation}
When $k=0$ and $\omega \neq 0$, we use
\begin{equation}
\bm{C}_j (A) 
 = \frac{\omega }{2\pi}\int^{2\pi /\omega }_{0}\! f(\sum^N_{n=1} A_n e^{ij_n \omega t}\bm{w}_n) e^{-ij\omega t} dt 
\end{equation}
instead of (\ref{4-47}).
Then, we obtain the expansion
\begin{equation}
f(\sum^N_{n=1} A_n e^{ij_n \omega t + ij_n kx}\bm{w}_n)
 = \sum^\infty_{j=-\infty} \bm{C}_j(A)e^{ij\omega t + ijkx}.
\end{equation}
Further, define the function $\bm{R}_n : (B^r)^N \to (B^r)^m$ to be
\begin{equation}
\bm{R}_n(A) = \Pi_n \bm{C}_{j_n} (A), \quad j_n\in J,\,\, (n=1,\cdots ,N),
\end{equation}
and define $R_n : (B^r)^N \to B^r$ so that $\bm{R}_n(A)=R_n(A)\bm{w}_n$.
Let $\varepsilon >0$ be a small parameter.
Let us consider the two initial value problems:
\begin{subnumcases}
{}
\displaystyle \frac{\partial u}{\partial t}  = \mathcal{P}u + \varepsilon f(u) , \quad
\displaystyle u(0,x) = \sum^N_{n=1} e^{ij_nkx} v_n(\eta \hat{x}_1)\bm{w}_n, \label{4-45a}\\[0.2cm]
\displaystyle \frac{\partial A_n}{\partial t} = \mathcal{Q}_nA_n + \varepsilon R_n(A), \quad
\displaystyle A_n(0,x) = v_n(\eta \hat{x}_1),\,\, (n=1,\cdots ,N),\label{4-45b}
\end{subnumcases}
where $\eta = \varepsilon ^{1/M}$ and $M :=\min \{ M_1,\cdots ,M_N\}$.
\\[0.2cm]
\textbf{Example \thedef.} 
Let us consider the system (\ref{1-14}), whose perturbation term is given by
\begin{equation}
f(u,v) = \left(
\begin{array}{@{\,}c@{\,}}
u-u^3 \\
0
\end{array}
\right) .
\end{equation}
The reduction of the linear part was calculated in Example 3.5, in which it was shown that
(E0) to (E3) are satisfied with $m=d=2,\, D=1,\, k=(0,c), \, \omega =0$,
$J = \{ j_1 = 1,\, j_2 = -1\}$ and $M_1 = M_2 = 4$.
We use the same notation as Example 3.5.
The matrix $P(i\xi)$ at $(\xi_1, \xi_2) = (0, \pm c)$ is given by
\begin{eqnarray*}
P(\pm i c) = \left(
\begin{array}{@{\,}cc@{\,}}
(k+d)/(2d) &  -1 \\
1 & (k+d)/2
\end{array}
\right).
\end{eqnarray*}
Eigenvalues and eigenvectors of this matrix are
\begin{eqnarray*}
& & \bm{w} = \left(
\begin{array}{@{\,}c@{\,}}
1 \\
(k+d)/2
\end{array}
\right) \quad \mathrm{for}\,\, \lambda _1 (0, \pm c) = 0, \\
& & \bm{v} = \left(
\begin{array}{@{\,}c@{\,}}
1 \\
(k+d)/(2d)
\end{array}
\right) \quad \mathrm{for}\,\, \lambda _2 (0, \pm c) = -\frac{(1-d)(k+d)}{2d} < 0.
\end{eqnarray*}
Since both eigenvectors of $\lambda _1(0, c)$ and $\lambda _1(0,-c)$ are given by $\bm{w}$ above,
we calculate the Fourier expansion of $f(A_1 e^{icy} \bm{w} + A_2 e^{-icy} \bm{w})$.
Then, it turns out that 
\begin{eqnarray*}
\bm{C}_1 (A) = \left(
\begin{array}{@{\,}c@{\,}}
A_1 - 3A_1^2A_2 \\
0
\end{array}
\right), \quad \bm{C}_{-1} (A) = \left(
\begin{array}{@{\,}c@{\,}}
A_2 - 3A_1A_2^2 \\
0
\end{array}
\right).
\end{eqnarray*}
Then, it is easy to show that projections of them are
\begin{eqnarray*}
\bm{R}_1(A) = \Pi_1 \bm{C}_1(A) = \frac{A_1 - 3A_1^2A_2}{1-d} \bm{w}, \quad
\bm{R}_2(A) = \Pi_2 \bm{C}_{-1}(A) = \frac{A_2 - 3A_1A_2^2}{1-d}\bm{w}.
\end{eqnarray*}
Therefore, the amplitude equation (\ref{4-45b}) is given by
\begin{equation}
\left\{ \begin{array}{l}
\displaystyle \frac{\partial A_1}{\partial t} = \mathcal{Q}A_1+\frac{\varepsilon }{1-d}(A_1 - 3A_1^2 A_2),  \\[0.2cm]
\displaystyle \frac{\partial A_2}{\partial t} = \mathcal{Q}A_2+\frac{\varepsilon }{1-d}(A_2 - 3A_1 A_2^2),  \\
\end{array} \right.
\end{equation}
where $\mathcal{Q}$ is defined by (\ref{3-17}).
If we suppose $A_2 = \overline{A}_1$, the equation (\ref{1-17}) is obtained.
\\[0.2cm]
\textbf{Example \thedef.} 
Noting Remark 3.7, to extend (E1) to the case that several eigenvalues lie on the imaginary axis is straightforward 
and the amplitude equation is defined in a similar manner as above.
Let us consider the system (\ref{1-7b}), whose perturbation term is the same as the above Example.
By the same calculation as Example 3.8 and 4.11, we obtain the amplitude equation 
\begin{equation}
\left\{ \begin{array}{l}
\displaystyle \frac{\partial A_1}{\partial t} = \mathcal{Q}A_1+\frac{\varepsilon }{2}(A_1 - 3A_1^2 A_2),  \\[0.2cm]
\displaystyle \frac{\partial A_2}{\partial t} = \mathcal{Q}A_2+\frac{\varepsilon }{2}(A_2 - 3A_1 A_2^2),  \\
\end{array} \right.
\end{equation}
where $\mathcal{Q} = D \partial ^2$ is obtained in Example 3.8.
If we suppose $A_2 = \overline{A}_1$, the equation (\ref{1-7c}) is obtained.
\\[0.2cm]
\textbf{Theorem \thedef.} 
Suppose (E0) to (E3), $f : (B^r)^m \to (B^r)^m\,\, (r\geq 1)$ is $C^1$ and $\varepsilon >0$ is sufficiently small.
For any $\{ v_n\}^N_{n=1} \subset B^{r}$,
there exist positive numbers $C,T_0$ and $t_0$ such that mild solutions of 
the two initial value problems (4.50) satisfy
\begin{equation}
|| u(t,x) - \sum^N_{n=1}A_n(t,x) e^{ij_n\omega t + ij_nkx} \bm{w}_n|| 
\leq C\eta = C\varepsilon ^{1/M},
\end{equation}
for $t_0 \leq t \leq T_0/\varepsilon $.
\\

Further, suppose that
\\
\textbf{(E4)} For $|\beta| = M_n$ and $ n=1,\cdots ,N$,
\[ P(i\xi) = \overline{P(-i\xi)},\,\, f(\overline{u}) = \overline{f(u)}\quad 
\text{and}  \quad \displaystyle \frac{\partial^\beta  \lambda _1}{\partial \hat{x}_1^\beta}(j_nk) = 
\frac{\partial^\beta  \lambda _1}{\partial \hat{x}_1^\beta}(-j_nk).\]

In this case, the set $J$ consists of $J = \{ j_1, \cdots ,j_N\} \cup \{ -j_1, \cdots , -j_N\}$ as is Sec.4.1.
We consider the two initial value problems:
\begin{subnumcases}
{}
\displaystyle \frac{\partial u}{\partial t}  = \mathcal{P}u + \varepsilon f(u) , \quad
\displaystyle u(0,x) = \sum^N_{n=1} \left( e^{ij_n kx}
                     +  e^{-ij_n kx}\right) v_n(\eta \hat{x}_1) \bm{w}_n, \label{4-52a}\\[0.2cm]
\displaystyle \frac{\partial A_n}{\partial t} = \mathcal{Q}_nA_n + \varepsilon S_n(A), \quad
\displaystyle A_n(0,x) = v_n(\eta \hat{x}_1),\,\, (n= 1,\cdots ,N), \label{4-52b}
\end{subnumcases}
where $S_n(A)$ is defined by (\ref{4-13b}).
In the next proposition, $B^r = BC^r(\mathbf{R}^d; \mathbf{R})$ denotes the set of real-valued functions.
\\[0.2cm]
\textbf{Theorem \thedef.}  
Suppose (E0) to (E4) and 
$f : (B^r)^m \to (B^r)^m\,\, (r\geq 1)$ is $C^2$ such that the second derivatives are locally Lipschitz continuous.
Suppose that there exists a constant vector $\phi = (\phi_1, \cdots ,\phi_N) \in \mathbf{R}^N$ such that
\\
\textbf{(i)} $S_n(\phi) = 0$ for $n=1,\cdots ,N$,
\\
\textbf{(ii)} the Jacobi matrix of $(S_1 ,\cdots ,S_N)$ at $\phi$ is diagonalizable and 
all eigenvalues of the matrix have negative real parts.
\\
If $\varepsilon >0$ is sufficiently small, Eq.(\ref{4-52a}) has a solution of the form
\begin{equation}
u_p(t,x) = \sum^N_{n=1} \Bigl( \phi_n \bm{w}_n + \eta \bm{\psi}_n (t,x, \eta) \Bigr) \cdot 
 (e^{ij_n\omega t + ij_nkx} + e^{-ij_n\omega t - ij_nkx}).
\end{equation}
The vector-valued functions $\bm{\psi}_n $ and $u_p$ are bounded as $\eta \to 0$ and satisfy
\begin{eqnarray*}
\left\{ \begin{array}{l}
\text{$2\pi / \omega $-periodic in $t$ (when $\omega \neq 0$)}, \\
\text{constant in $t$ (when $\omega = 0$)},  \\
\end{array} \right.
\quad
\left\{ \begin{array}{l}
\text{$2\pi / k_j $-periodic in $x_j$ (when $k_j \neq 0$)}, \\
\text{constant in $x_j$ (when $k_j = 0$)},  \\
\end{array} \right.
\end{eqnarray*}
for $j=1,\cdots ,d$.
This $u_p$ is stable in the following sense:
For any $n=1,\cdots ,N$, there is a neighborhood $U_n\subset B^r$ of $\phi_n$ in $B^r$ such that if $v_n\in U_n$,
then a mild solution $u$ of the initial value problem (\ref{4-52a}) satisfies
$|| u(t, \cdot) - u_p(t, \cdot) || \to 0$ as $t\to \infty$.
\\[0.2cm]
\textbf{Proof.}
We suppose $D=d$ for simplicity and use the same notation as the proof of Thm.4.2 (see Step 1).
Two mild solutions of Eq.(4.50) satisfy
\begin{eqnarray*}
& & \hat{u} - \sum^N_{n=1} \hat{A}_n e^{ij_n\omega T/\varepsilon + ij_nkX/\eta }\bm{w}_n 
= \sum^N_{n=1}e^{\hat{\mathcal{P}}T} (e^{ij_nkX/\eta }v_n \bm{w}_n) 
         -\sum^N_{n=1}e^{ij_n\omega T/\varepsilon + ij_nkX/\eta }(e^{\hat{\mathcal{Q}}_nT}v_n) \bm{w}_n  \\
& &  \quad\quad        + \int^{T}_{0}\! e^{\hat{\mathcal{P}}(T-s)} f(\hat{u}(s))ds
       - \sum^N_{n=1} \int^{T}_{0}\!e^{ij_n\omega T/\varepsilon + ij_nkX/\eta } 
                   e^{\hat{\mathcal{Q}}_n(T-s)}R_n(\hat{A}(s)) \bm{w}_n ds \\
&=& F(T) + \int^{T}_{0}\! e^{\hat{\mathcal{P}}(T-s)} 
\Bigl(f(\hat{u}(s))-f(\textstyle \sum^N_{n=1} \hat{A}(s)e^{ij_n\omega s/\varepsilon  + ij_nkX/\eta }\bm{w}_n)\Bigr) ds,
\end{eqnarray*}
where 
\begin{eqnarray*}
& & F(T) = \sum^N_{n=1}e^{\hat{\mathcal{P}}T} (e^{ij_nkX/\eta }v_n \bm{w}_n) 
         -\sum^N_{n=1}e^{ij_n\omega T/\varepsilon + ij_nkX/\eta }(e^{\hat{\mathcal{Q}}_nT}v_n) \bm{w}_n \\
&+& \!\!\! \int^{T}_{0}\!\! e^{\hat{\mathcal{P}}(T-s)}
f(\textstyle \sum^N_{n=1} \hat{A}(s)e^{ij_n\omega s/\varepsilon  + ij_nkX/\eta } \bm{w}_n) ds \displaystyle
-\sum^N_{n=1} \int^{T}_{0}\!\!
     e^{ij_n\omega T/\varepsilon + ij_nkX/\eta } e^{\hat{\mathcal{Q}}_n(T-s)}R_n(\hat{A}(s)) \bm{w}_n ds.
\end{eqnarray*}
As before, Gronwall inequality yields
\begin{eqnarray*}
|| \hat{u} - \sum^N_{n=1} \hat{A}_n e^{ij_n\omega T/\varepsilon + ij_nkX/\eta } \bm{w}_n ||_\eta \leq 
|| F(T) ||_\eta + D_1L \int^{T}_{0}\! e^{D_1L(T-s)} || F(s)||_\eta ds
\end{eqnarray*}
for some constants $D_1, L > 0$.
By using definitions of $\bm{C}_j$ and $R_n$, we rewrite $F(T)$ as
\begin{eqnarray}
F(T) &=& \sum^N_{n=1}e^{\hat{\mathcal{P}}T} (e^{ij_nkX/\eta }v_n \bm{w}_n) 
         -\sum^N_{n=1}e^{ij_n\omega T/\varepsilon + ij_nkX/\eta }(e^{\hat{\mathcal{Q}}_nT}v_n)\bm{w}_n \nonumber \\
&+& \sum^N_{n=1} \int^{T}_{0}\! e^{ij_n \omega s/\varepsilon } 
    \left( e^{\hat{\mathcal{P}}(T-s)} e^{ij_n kX/\eta} - e^{ij_n\omega (T-s)/\varepsilon + ij_nkX/\eta  } 
               e^{\hat{\mathcal{Q}}_n(T-s)}\right) R_n(\hat{A}(s)) \bm{w}_n ds \nonumber \\
&+& \sum_{j\notin J} H_j(T,X) \nonumber \\
&+& \sum^N_{n=1} \int^{T}_{0}\! e^{ij_n \omega s/\varepsilon } e^{\hat{\mathcal{P}}(T-s)}
        e^{ij_n kX/\eta} \Pi_n^{\bot} \bm{C}_{j_n} (\hat{A}(s))ds,
\label{4-54}
\end{eqnarray}
where $H_j$ is defined by
\begin{eqnarray}
& & H_j(T,X) = \int^{T}_{0}\!  e^{\hat{\mathcal{P}}(T-s)} e^{ij\omega s/\varepsilon + ijkX/\eta}
                 \bm{C}_j(\hat{A}(s))ds.
\end{eqnarray}
Now we have arrived at the same situation as (\ref{4-20b}) except for the last term.
\\[0.2cm]
\textbf{Lemma \thedef.} 
There exists a function $g_n : [0, T_0] \times \mathbf{R}^d \times C([0,T_0];B^r)^N \to C([0,T_0];B^r)^m$ such that
\begin{equation}
\int^{T}_{0}\! e^{ij_n \omega s/\varepsilon } e^{\hat{\mathcal{P}}(T-s)}
        e^{ij_n kX/\eta} \Pi_n^{\bot} \bm{C}_{j_n} (\hat{A}(s))ds = \eta g_n(T,X, \hat{A}).
\end{equation}
$g(T,X, \hat{A})$ is Lipschitz continuous in $\hat{A} \in C([0,T_0];B^r)^N$ and bounded as $\eta \to 0$.
\\

If this lemma is true, the rest of the proofs of Thm.4.13, 4.14 are completely the same as 
those of Thm.4.2, 4.3.
\\[0.2cm]
\textbf{Proof.} 
We calculate $e^{\hat{\mathcal{P}}(T-s)}e^{ij_n kX/\eta} \Pi_n^{\bot} \bm{C}_{j_n} (\hat{A}(s))$ as
\begin{eqnarray*}
e^{\hat{\mathcal{P}}(T-s)}e^{ij_n kX/\eta} \Pi_n^{\bot} \bm{C}_{j_n} (\hat{A}(s))
&=& \!\!\frac{1}{(2\pi)^d} \!\int \!\!\! \int \! e^{-iy\xi} e^{P(i\eta \xi) (T-s)/\varepsilon }
 e^{ij_n k(X+y)/\eta }\Pi_n^{\bot} \bm{C}_{j_n} (\hat{A}(s, X+y))d\xi dy \\
&=& \!\frac{e^{ij_n kX/\eta }}{(2\pi)^d} \int \!\!\! \int \! e^{-iy\xi} e^{P(i\eta \xi + ij_nk) (T-s)/\varepsilon }
 \Pi_n^{\bot} \bm{C}_{j_n} (\hat{A}(s, X+y))d\xi dy.
\end{eqnarray*}
Let us suppose that the coordinate of $u = (u_1 , \cdots ,u_m)$ is defined so that the $m\times m$ matrix
$P(ij_nk)$ is diagonal.
Define a matrix $S(\xi)$ such that Eq.(\ref{3-22}) holds.
By the assumption, $S(j_nk)$ is the identity matrix.
Then, we obtain
\begin{eqnarray*}
& & e^{\hat{\mathcal{P}}(T-s)}e^{ij_n kX/\eta} \Pi_n^{\bot} \bm{C}_{j_n} (\hat{A}(s)) \\
&=& \frac{e^{ij_n kX/\eta }}{(2\pi)^d} \int \!\!\! \int \! 
            e^{-iy\xi}S(j_nk + \eta \xi) e^{\Lambda (j_nk + \eta \xi) (T-s)/\varepsilon }S(j_nk + \eta \xi)^{-1}
            \Pi_n^{\bot} \bm{C}_{j_n} (\hat{A}(s, X+y))d\xi dy.
\end{eqnarray*}
By expanding $S(j_nk + \eta \xi)^{\pm 1}$, it turns out that there is a function $G_n(T,X,\hat{A})$,
which is Lipschitz continuous in $\hat{A}\in C([0,T_0];B^r)^N$ and bounded as $\eta \to 0$, such that
\begin{eqnarray*}
& & e^{\hat{\mathcal{P}}(T-s)}e^{ij_n kX/\eta} \Pi_n^{\bot} \bm{C}_{j_n} (\hat{A}(s)) \\
&=& \frac{e^{ij_n kX/\eta }}{(2\pi)^d} \int \!\!\! \int \! 
       e^{-iy\xi} e^{\Lambda (j_nk + \eta \xi) (T-s)/\varepsilon } \Pi_n^{\bot} \bm{C}_{j_n} (\hat{A}(s, X+y))d\xi dy
 + \eta G_n(T-s,X,\hat{A}).
\end{eqnarray*}
Let $\bm{C}_{j_n}^{(l)}$ be the $l$-th component of the vector $\bm{C}_{j_n}$.
Due to the definition of $\Pi^{\bot}_n$, the first component of 
$e^{\Lambda (j_nk + \eta \xi) (T-s)/\varepsilon } \Pi_n^{\bot} \bm{C}_{j_n}$ is zero,
and the $l$-th component is given by
$e^{\lambda _l(j_nk + \eta \xi) (T-s)/\varepsilon } \bm{C}^{(l)}_{j_n}$ for $l=2,\cdots ,m$.
To prove the lemma, it is sufficient to estimate
\begin{eqnarray*}
I_{n,l} := \int^{T}_{0}\! e^{ij_n \omega s/\varepsilon } 
 \int \!\!\! \int \! 
     e^{-iy\xi} e^{\lambda_l (j_nk + \eta \xi) (T-s)/\varepsilon } \bm{C}^{(l)}_{j_n} (\hat{A}(s, X+y))d\xi dy ds.
\end{eqnarray*}
Because of (E1) and (E3), there exists a positive number $\beta$ such that
\begin{equation}
\mathrm{Re}[\lambda_l (j_nk + \eta \xi) + \beta ] \leq 0, \quad 
\mathrm{Re}[\lambda_l (j_nk + \eta \xi) + \beta ] \sim O(-|\xi|^2)
\end{equation}
as $|\xi| \to \infty$ for any $l= 2,\cdots ,m$.
Then, we obtain
\begin{eqnarray*}
I_{n,l} &=& \int^{T}_{0}\! e^{ij_n \omega s/\varepsilon } e^{-\beta (T-s)/\varepsilon }
 \int \!\!\! \int \! e^{-iy\xi} 
    e^{(\lambda_l (j_nk + \eta \xi)+\beta ) (T-s)/\varepsilon } \bm{C}^{(l)}_{j_n} (\hat{A}(s, X+y))d\xi dy ds \\
&:=& \int^{T}_{0}\!  e^{ij_n \omega s/\varepsilon } e^{-\beta (T-s)/\varepsilon }K_{n,l}(T-s,X,\hat{A})ds,
\end{eqnarray*}
where $K_{n,l}$ is Lipschitz continuous in $\hat{A} \in C([0,T_0]; B^r)^N$ and bounded as $\eta \to 0$.
The mean value theorem proves that there exists $0\leq \tau \leq T$ such that
\begin{eqnarray*}
I_{n,l} =  e^{ij_n \omega \tau /\varepsilon }K_{n,l}(T-\tau,X,\hat{A}(\tau, X)) 
\int^{T}_{0}\!  e^{-\beta (T-s)/\varepsilon }ds,
\end{eqnarray*}
which is of order $O(\varepsilon )$.
Hence, putting $(0, I_{n,2}, \cdots ,I_{n,N}) = \varepsilon I_n(T,X,\hat{A})$ and 
\begin{eqnarray*}
\eta g_n(T,X, \hat{A}) = \varepsilon  \frac{e^{ij_nkX/\eta}}{(2\pi)^d}I_n(T,X,\hat{A})
  + \eta \int^{T}_{0}\! e^{ij_n \omega s /\varepsilon }G_n(T-s, X, \hat{A})ds 
\end{eqnarray*}
proves the lemma. \hfill $\blacksquare$
\\

Now the function $F(T)$ in Eq.(\ref{4-54}) is estimated with the aid of Prop.3.6, Lemma 4.4 and Lemma 4.15
to show $|| F(T) ||_{\eta} \sim O(\eta)$.
Then, the Gronwall inequality proves Thm.4.13.
A proof of Thm.4.14 is also done in the same way as that of Thm.4.3. \hfill $\blacksquare$

%%%%%%%%%%%%%%%%%%%%%%%%%%%%%%%%%%%%%%%%%%%%%%%%%%%%%%%%%%%%%%%%%%%%%%%%%%%%%%%%%%%%%%%%%%%%%%%%%%%%%%%%%%%%
%%%%%%%%%%%%%%%%%%%%%%%%%%%%%%%%%%%%%%%%%%%%%%%%%%%%%%%%%%%%%%%%%%%%%%%%%%%%%%%%%%%%%%%%%%%%%%%%%%%%%%%%%%%%

\vspace*{0.5cm}
\textbf{Acknowledgements.}

The author would like to thank Professor Yoshihisa Motira and Shin-Ichiro Ei for useful comments.
This work was supported by Grant-in-Aid for Young Scientists (B), No.22740069 from MEXT Japan.

%%%%%%%%%%%%%%%%%%%%%%%%%%%%%%%%%%%%%%%%%%%%%%%%%%%%%%%%%%%%%%%%%%%%%%%%%%%%%%%%%%%%%%%%%%%%%%%%%%
%%%%%%%%%%%%%%%%%%%%%%%%%%%%%%%%%%%%%%%%%%%%%%%%%%%%%%%%%%%%%%%%%%%%%%%%%%%%%%%%%%%%%%%%%%%%%%%%%%

\end{document}